\documentclass{icmart}

\contact[j.marklof@bristol.ac.uk]{School of Mathematics, University of Bristol, Bristol BS8 1TW, U.K.}

\usepackage[active]{srcltx}

\usepackage{graphicx}
\usepackage{hyperref}
\usepackage{enumerate}
\usepackage{bbm}

\newtheorem{thm}{Theorem}
\theoremstyle{remark}
\theoremstyle{plain}

\numberwithin{equation}{section}

\def\NN{{\mathbb N}}
\def\QQ{{\mathbb Q}}
\def\PP{{\mathbb P}}
\def\RR{{\mathbb R}}
\def\TT{{\mathbb T}}

\def\ZZ{{\mathbb Z}}

\def\mm{{\mathbbm m}}
\def\pp{{\mathbbm p}}

\def\BG{{\operatorname{BG}}}

\def\veca{{\text{\boldmath$a$}}}
\def\vecb{{\text{\boldmath$b$}}}

\def\vece{{\text{\boldmath$e$}}}
\def\vech{{\text{\boldmath$h$}}}

\def\vecell{{\text{\boldmath$\ell$}}}

\def\vecq{{\text{\boldmath$q$}}}
\def\vecQ{{\text{\boldmath$Q$}}}

\def\vecs{{\text{\boldmath$s$}}}

\def\vect{{\text{\boldmath$t$}}}

\def\vecv{{\text{\boldmath$v$}}}
\def\vecV{{\text{\boldmath$V$}}}
\def\vecw{{\text{\boldmath$w$}}}

\def\vecx{{\text{\boldmath$x$}}}

\def\vecy{{\text{\boldmath$y$}}}

\def\vecalf{{\text{\boldmath$\alpha$}}}

\def\ueta{{\underline{\eta}}}
\def\uxi{{\underline{\xi}}}

\def\vecnull{{\text{\boldmath$0$}}}

\def\scrA{{\mathcal A}}
\def\scrB{{\mathcal B}}
\def\scrC{{\mathcal C}}
\def\scrD{{\mathcal D}}

\def\scrF{{\mathcal F}}

\def\scrK{{\mathcal K}}
\def\scrL{{\mathcal L}}

\def\scrP{{\mathcal P}}

\def\scrT{{\mathcal T}}

\def\scrV{{\mathcal V}}
\def\scrW{{\mathcal W}}

\def\fO{{\mathfrak O}}

\def\fZ{{\mathfrak Z}}

\def\e{\mathrm{e}}

\def\id{\operatorname{id}}
\def\intl{{\operatorname{int}}}

\def\C{\operatorname{C{}}}

\def\L{\operatorname{L{}}}
\def\M{\operatorname{M{}}}

\def\S{\operatorname{S{}}}

\def\SL{\operatorname{SL}}

\def\SO{\operatorname{SO}}

\def\T{\operatorname{T{}}}

\def\vol{\operatorname{vol}}

\def\GamG{\Gamma\backslash G}

\def\SLdZ{\SL(d,\ZZ)}
\def\SLdR{\SL(d,\RR)}

\def\trans{\,^\mathrm{t}\!}

\def\nbar{\overline{n}}
\def\xibar{\overline{\xi}}

\def\hatw{{\widehat{\vecw}}}

\def\bs{\backslash}

\def\UB{{\scrB_1^{d-1}}}
\def\US{{\S_1^{d-1}}}

\def\Exp{\operatorname{\mathbf{E}}}

\title[The low-density limit of the Lorentz gas]{The low-density limit of the Lorentz gas: periodic, aperiodic and random}

\author[Jens Marklof]{Jens Marklof\thanks{The research leading to the results presented here has received funding from the European Research Council under the European Union's Seventh Framework Programme (FP/2007-2013) / ERC Grant Agreement n. 291147. The author is furthermore supported by a Royal Society Wolfson Research Merit Award.}}

\begin{document}

\begin{abstract}
The Lorentz gas is one of the simplest, most widely used models to study the transport properties of rarified gases in matter. It describes the dynamics of a cloud of non-interacting point particles in an infinite array of fixed spherical scatterers. More than one hundred years after its conception, it is still a major challenge to understand the nature of the kinetic transport equation that governs the macroscopic particle dynamics in the limit of low scatterer density (the Boltzmann-Grad limit). Lorentz suggested that this equation should be the linear Boltzmann equation. This was confirmed in three celebrated papers by Gallavotti, Spohn, and Boldrighini, Bunimovich and Sinai, under the assumption that the distribution of scatterers is sufficiently disordered. In the case of strongly correlated scatterer configurations (such as crystals or quasicrystals), we now understand why the linear Boltzmann equation fails and what to substitute it with. A particularly striking feature of the periodic Lorentz gas is a heavy tail for the distribution of free path lengths, with a diverging second moment, and superdiffusive transport in the limit of large times. 
\end{abstract}



\maketitle


\section{Introduction}\label{sec:introduction}

The Lorentz gas describes the time evolution of a cloud of non-interacting point particles in an infinitely extended array of fixed scatterers. In the simplest setting of zero external force fields, each particle moves with constant velocity along a straight line until it hits a sphere of radius $r$, where it is scattered elastically. Besides specular reflection (as in Lorentz' original setting), we will also allow more general spherically symmetric scattering maps, for example those resulting from muffin-tin Coulomb potentials.  
The scatterers are centered at the points of a locally finite subset $\scrP\subset\RR^d$, which is fixed once and for all. The configuration space of the Lorentz gas is thus $\scrK_r= \RR^d \setminus (\scrP+\scrB_r^d)$
where $\scrB_r^d$ is the open ball in $\RR^d$ of radius $r$, centered at the origin. The phase space of the Lorentz gas is $\T(\scrK_r)$, the tangent bundle of $\scrK_r$. We use the convention that, for $\vecq\in\partial\scrK_r$, the tangent vector $\vecv$ points away from the scatterer.\footnote{We ignore the case when scatterers overlap. This configuration will be statistically insignificant in the limit $r\to 0$ for $\scrP$ with constant density.} Given initial data $(\vecq,\vecv)\in \T(\scrK_r)$ at time $t=0$, we denote position and velocity at time $t\in\RR$ by $(\vecq(t),\vecv(t))$. For notational reasons it is convenient to also define the dynamics inside the scatterer by 
$(\vecq(t),\vecv(t))=(\vecq,\vecv)$ for every $(\vecq,\vecv)\in \T(\RR^d)\setminus\T(\scrK_r)$. With this, the phase space is $\T(\RR^d)=\RR^d\times\RR^d$. The Liouville measure of our dynamics is the Lebesgue measure $d\vecq\,d\vecv$. Since we have assumed that the scattering map is elastic, the particle speed $\|\vecv\|$ is a constant of motion. We may therefore restrict the dynamics, without loss of generality,  to the unit tangent bundle $\T^1(\RR^d)=\RR^d\times\S_1^{d-1}$, where the Liouville measure is now the Lebesgue measure restricted to $\|\vecv\|=1$. 

We assume that $\scrP$ has constant density $\nbar>0$, i.e.\ for any bounded $\scrD\subset\RR^d$ with $\vol_{\RR^d}(\scrD)>0$ and $\vol_{\RR^d}(\partial\scrD)=0$ ($\vol_{\RR^d}$ denotes the Lebesgue measure in $\RR^d$ and $\partial\scrD$ the boundary of $\scrD$) we have
\begin{equation}\label{density}
\lim_{R\to\infty} \frac{\#(\scrP\cap R\scrD)}{\vol_{\RR^d}(R\scrD)} =\nbar.
\end{equation}
By a trivial rescaling of length units, we may assume in the following that $\nbar=1$.
\begin{figure}
\begin{center}
\begin{minipage}{0.95\textwidth}
\unitlength0.1\textwidth
\begin{picture}(10,6)(0,0)
\put(0,0){\includegraphics[width=\textwidth]{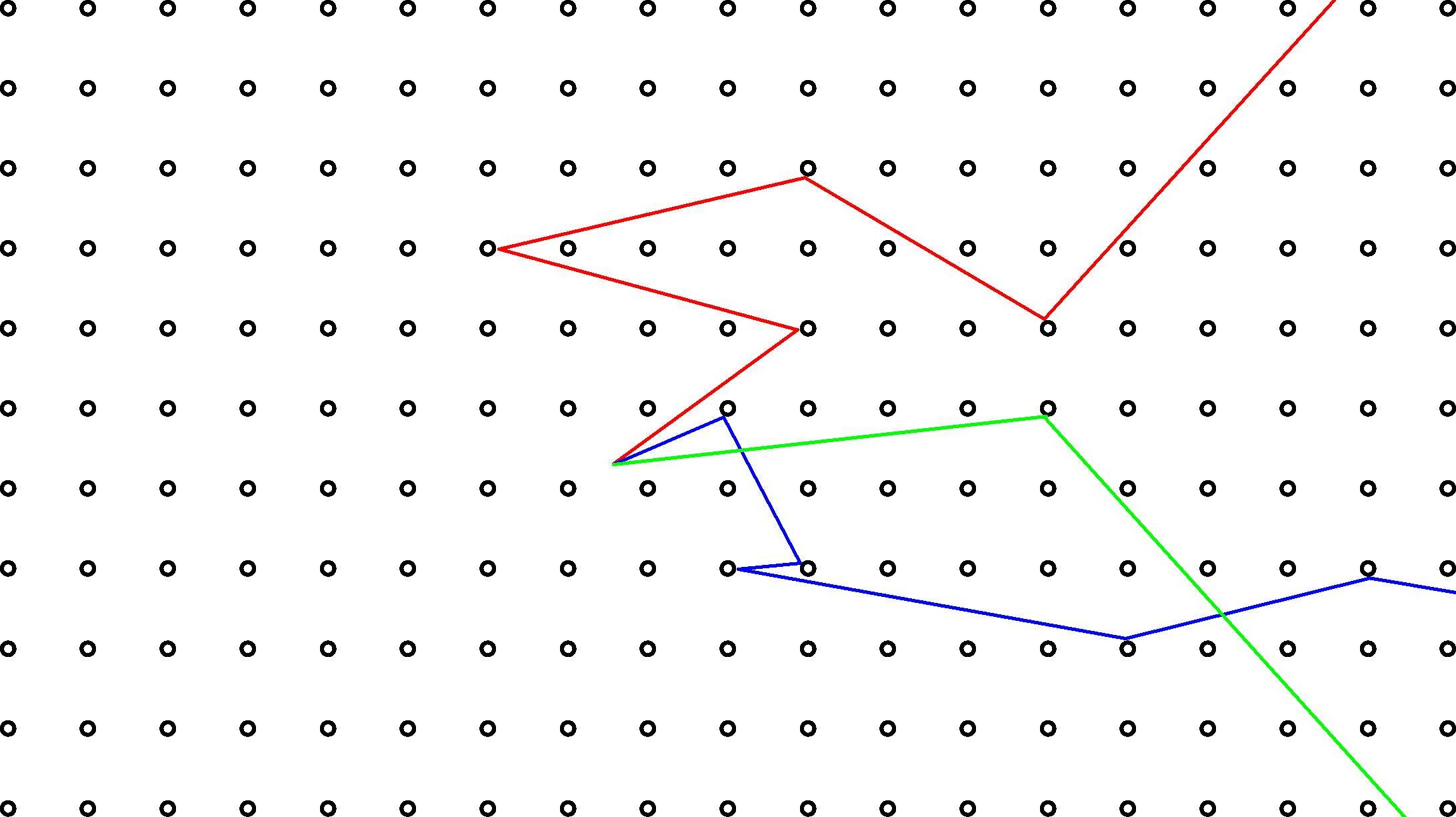}}
\put(3.9,2.4){$\vecq_0$}
\end{picture}
\end{minipage}
\end{center}
\caption{The Lorentz gas for a periodic scatterer configuration $\scrP=\ZZ^2$, with three distinct particle trajectories starting at the point $\vecq_0$. \label{figPL}}
\end{figure}

In the present setting, the Boltzmann-Grad limit is defined as the limit of low scatterer density. {\em Density} refers here to the {\em volume} density, i.e., the relative volume $v_d r^d$ occupied by the scatterers, rather than their {\em number} density $\nbar=1$. The constant $v_d=\vol_{\RR^d}(\scrB_1^d)=\pi^{d/2}/\Gamma(\tfrac{d+2}{2})$ is the volume of the $d$-dimensional unit ball. 
For a fixed scatterer configuration $\scrP$ the Boltzmann-Grad limit corresponds therefore to taking $r\to 0$. To capture the dynamics of the Lorentz gas in this limit, we measure length and time in units of the mean free path length,\footnote{The mean free path length is defined as the average distance travelled between collisions.} which is asymptotic to $v_{d-1}^{-1} r^{d-1}$ (as $r\to 0$). 
To this end we introduce the {\em macroscopic} coordinates 
\begin{equation}\label{macrocoo}
(\vecQ(t),\vecV(t)) = (r^{d-1} \vecq(r^{-(d-1)} t), \vecv(r^{-(d-1)} t) ) \in\T^1(\RR^d).
\end{equation}
The mean free path length is now given by the $r$-independent quantity $\xibar=v_{d-1}^{-1}$.
The evolution of an initial macroscopic particle density $f\in\T^1(\RR^d)$ is defined by the linear operator 
\begin{equation}
[L_r^t f](\vecQ,\vecV) := f(\vecQ(-t),\vecV(-t))
\end{equation}
where $(\vecQ(-t),\vecV(-t))$ are the macroscopic particle coordinates corresponding to the data $(\vecQ(0),\vecV(0))=(\vecQ,\vecV)$ at time $t=0$. 

The question is: For a given scatterer configuration $\scrP$, does $L_r^t$ have a (weak) limit as $r\to 0$? That is, for every $t>0$ is there 
\begin{equation}
L^t:\L^1(\T^1(\RR^d))\to\L^1(\T^1(\RR^d))
\end{equation}
such that, for every $f\in\L^1(\T^1(\RR^d))$ and bounded $\scrA\subset\T^1(\RR^d)$ with boundary of zero Lebesgue measure,
\begin{equation}\label{Lorentz-limit}
\lim_{r\to 0} \int_\scrA L_r^t f(\vecQ,\vecV) \, d\vecQ\,d\vecV = \int_\scrA L^t f(\vecQ,\vecV) \, d\vecQ\,d\vecV \; ?
\end{equation}
Using Boltzmann's heuristics, Lorentz \cite{Lorentz05} predicted in 1905 that the answer to this question should be ``yes'' and that the particle density $f_t:=L^t f$ at time $t$ satisfies the {\em linear} Boltzmann equation (also referred to as the {\em kinetic Lorentz equation})
\begin{equation}\label{LBeq}
(\partial_t +\vecV\cdot\partial_\vecQ) f_t(\vecQ,\vecV) 
=  \int_{\RR^d} [f_t(\vecQ,\vecV')-f_t(\vecQ,\vecV)] \,\sigma(\vecV,\vecV')\, d\vecV',
\end{equation}
where $\sigma(\vecV,\vecV')$ is the differential cross section of a single scatterer (see Section \ref{sec:intercollision}).
Lorentz' heuristic derivation was, over sixty years later, confirmed rigorously for random scatterer configurations $\scrP$ by Gallavotti \cite{Gallavotti69} and Spohn \cite{Spohn78}, where the convergence in \eqref{Lorentz-limit} is established for the ensemble average. Boldrighini, Bunimovich and Sinai \cite{Boldrighini83} proved a stronger result by showing that for a {\em fixed} realisation of a Poisson process the limit \eqref{Lorentz-limit} exists almost surely (cf.~Section \ref{sec:Poisson}). One can in fact show that, for initial data $(\vecQ_0,\vecV_0)$ randomly distributed in $\T^1(\RR^d)$ according to an absolutely continuous probability measure $\Lambda$, the curve $t\mapsto (\vecQ(t),\vecV(t))$ converges in distribution to a random flight process, where the free flight times are independent identically distributed random variables with an exponential distribution. Eq.~\eqref{LBeq} is precisely the Fokker-Planck-Kolmogorov equation of the limit process (cf.~Section~\ref{sec:Poisson}). 

In his 2006 ICM address \cite{Golse06}  (cf.~also \cite{Golse08}), Golse pointed out  that, due to the heavy tail of the free path length distribution \cite{Bourgain98,Golse00,partIV}, the linear Boltzmann equation fails in the case $\scrP=\ZZ^d$. The main objective of this paper is to illustrate the deeper reason behind this failure not only for general {\em periodic} scatterer configurations, see Section \ref{sec:periodic} and \cite{Caglioti10,partIII,partI,partII,partIV}, but as well for {\em aperiodic}  point sets with strong long-range correlations, cf.~Sections \ref{sec:union}, \ref{sec:quasicrystals} and \cite{union,quasi,Wennberg12}. We will uncover a new class of random flight processes that emerge in the Boltzmann-Grad limit (Sections \ref{sec:intercollision}, \ref{sec:beyond}) and whose transport equations generalise the linear Boltzmann equation \eqref{LBeq} in a natural way (Section \ref{sec:generalized}). 

A major open question in the field is whether the dynamics in the Lorentz gas converges, in the limit of large times $t$, to Brownian motion. The first seminal result in this direction was the proof of a central limit theorem for the two-dimensional {\em periodic} Lorentz gas with finite horizon\footnote{{\em Finite horizon} means that the free path length has an upper bound. This requires a suitable choice of scatterer configuration $\scrP$ (e.g.\ a triangular lattice) and sufficiently large scatterer radius~$r$.} by Bunimovich and Sinai \cite{Bunimovich:1980ur}. For general invariance principles, see Melbourne and Nicol \cite{Melbourne:2009ju} and references therein. In the case of the infinite-horizon periodic Lorentz gas, again in dimension $d=2$ and with fixed radius $r>0$, Bleher \cite{Bleher:1992ku} conjectured a superdiffusive central limit theorem with a $\sqrt{t \log t}$ normalisation, rather than the standard $\sqrt t$ in the finite horizon case. Bleher's conjecture was first proved by Sz\'asz and Varj\'u \cite{Szasz:2007uo} for the discrete-time billiard map, and by Dolgopyat and Chernov \cite{Dolgopyat:2009bl} for the billiard flow.\footnote{Superdiffusive central limit theorems have also been established for compact planar billiards, such as the stadium  \cite{Balint:2006jv} and billiards with cusps \cite{Balint:2011jt}.} It is currently unknown how to extend these results to higher dimensions $d\geq 3$ or to aperiodic scatterer configurations \cite{Balint:2008kt,Balint:2012fg,Chernov:1994hb,Dettmann12,Dettmann14,Lenci11,Nandori12,Szasz08}. The problem becomes tractable, however, if we pass to the low-density limit $r\to 0$: If $\scrP$ is a typical realisation of a Poisson process, then the limiting random flight process satisfies a central limit theorem with $\sqrt t$ scaling, in any dimension $d\geq 2$. This follows from standard techniques in the theory of Markov processes \cite{Papanicolaou75} as pointed out by Spohn \cite{Spohn78}. If $\scrP$ is a Euclidean lattice, then the limiting random flight process satisfies a {\em superdiffusive} central limit theorem with $\sqrt{t\log t}$ normalisation, again in any dimension $d\geq 2$. See Section \ref{sec:superdiffusion} and \cite{super} for further details.

\section{Intercollision flights}\label{sec:intercollision}

\begin{figure}
\begin{center}
\begin{minipage}{0.9\textwidth}
\unitlength0.1\textwidth
\begin{picture}(10,5)(0,0)
\put(0,0){\includegraphics[width=\textwidth]{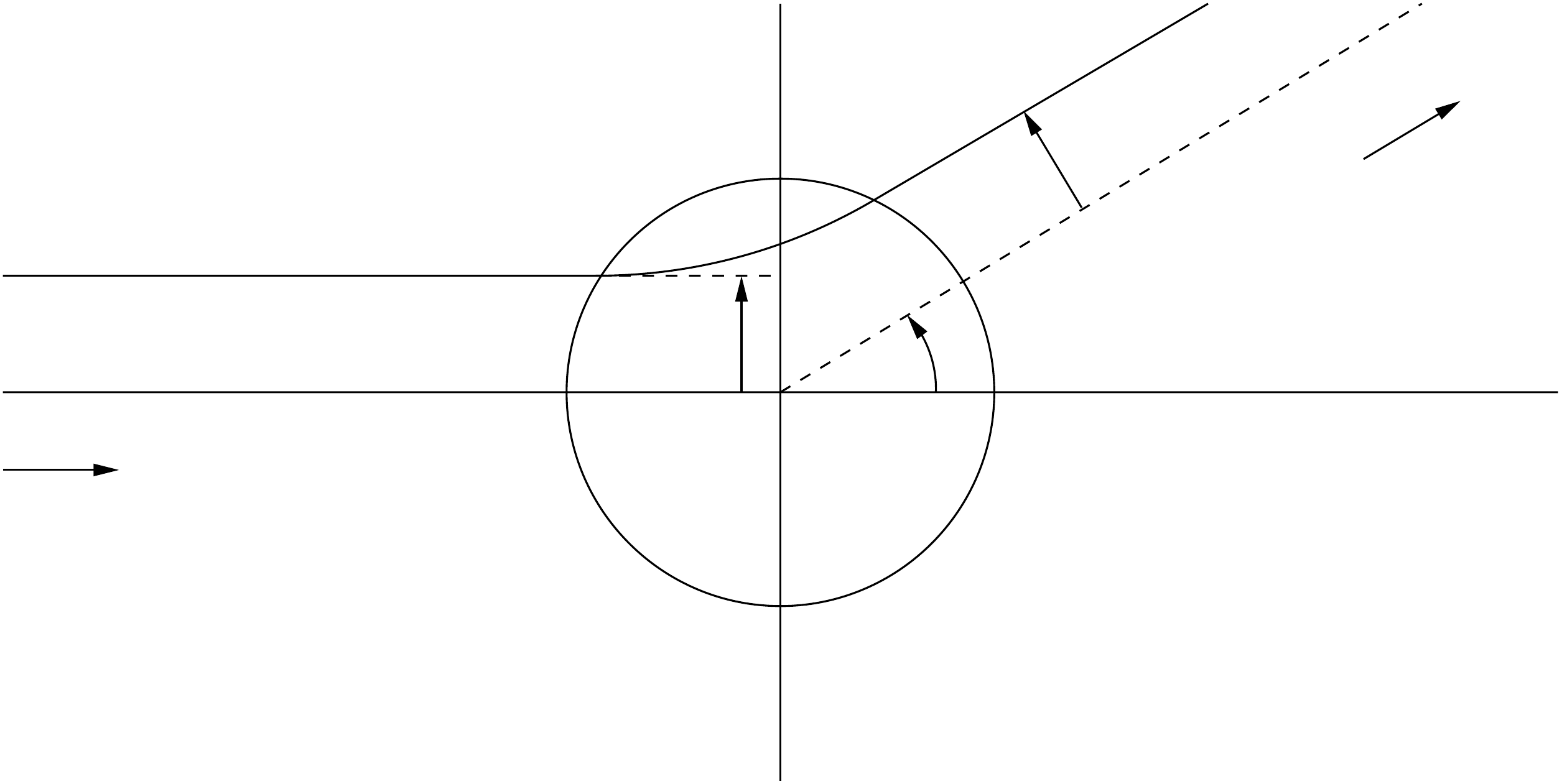}}
\put(0,1.75){$\vecv_\text{in}$} 
\put(4.33,2.75){$\vecb$}
\put(5.6,2.6){$\theta$} 
\put(6.85,4){$\vecs$}
\put(9.25,4){$\vecv_\text{out}$}
\end{picture}
\end{minipage}
\end{center}
\caption{Scattering in the unit ball.} \label{fig1}
\end{figure}
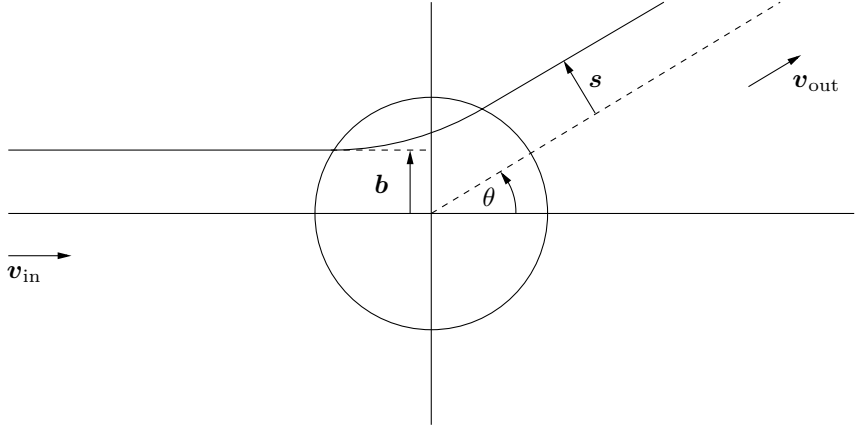

We begin by defining the scattering map, which we assume is spherically symmetric, preserves angular momentum and is the same for each scatterer. Let us choose a coordinate frame so that the incoming velocity is aligned with the first coordinate axis (cf.~Figure~\ref{fig1}),
\begin{equation}
\vecv_\text{in} =\vece_1:= (1, 0 , \ldots , 0 ).
\end{equation}
(All vectors are represented as row vectors.)
The {\em impact parameter} $\vecb$ is the orthogonal projection of the point of impact onto the plane orthogonal to $\vecv_\text{in}$, measured in units of $r$. In the present frame, $\vecb=(0,\vecw)$ with $\vecw\in\scrB_1^{d-1}$. (We will also refer to $\vecw$ as impact parameter.) When $\vecw\neq\vecnull$, the outgoing velocity is
\begin{equation}\label{inout}
\vecv_\text{out} =  \vecv_\text{in} \cos\theta +  (0,\hatw) \sin\theta  ,
\end{equation}
where the angle $\theta$ is called the {\em scattering angle} and $\hatw:=w^{-1}\vecw$ with $w:=\|\vecw\|$. For $\vecw=\vecnull$ we simply assume $\vecv_\text{out} =  -\vecv_\text{in}$. By the assumed spherical symmetry, $\theta=\theta(w)$ is only a function of the length $w\in[0,1[$ of the impact parameter $\vecw$. Equation \eqref{inout} can be expressed as
\begin{equation}\label{inout2}
\vecv_\text{out} =  \vecv_\text{in} S(\vecw)^{-1}  ,
\end{equation}
with the matrix
\begin{equation}\label{SbLG}
S(\vecw) = \exp\begin{pmatrix} 0 & -\theta(w) \hatw \\ \theta(w) \trans\hatw & 0_{d-1} \end{pmatrix}\in\SO(d).
\end{equation}
The {\em exit parameter} is defined as the orthogonal projection of the point of exit onto the plane orthogonal to $\vecv_\text{out}$, and is given by
\begin{equation}\label{inout2exit}
\vecs = -w \vecv_\text{in} \sin\theta +  (0,\vecw) \cos\theta = (0,\vecw) S(\vecw)^{-1}  .
\end{equation}
The differential scattering cross section $\sigma(\vecv_\text{in},\vecv_\text{out})$ is defined by the relation
\begin{equation}\label{dcs}
\sigma(\vecv_\text{in},\vecv_\text{out}) \,d\vecv_\text{out} = d\vecw.
\end{equation}
Note that in the present setting $\sigma(\vecv_\text{in},\vecv_\text{out})=\sigma(\vecv_\text{out},\vecv_\text{in})$.

For simplicity, we assume throughout this paper that one of the following conditions holds:\footnote{All results extend in fact to more general scattering maps, see \cite{partII} for details.}

\begin{itemize}
\item[(A)] $\theta\in\C^1([0,1[\,)$ is strictly decreasing with $\theta(0)=\pi$ and $\theta(w)> 0$. 
\item[(B)] $\theta\in\C^1([0,1[\,)$ is strictly increasing with $\theta(0)=-\pi$  and $\theta(w)< 0$. 
\end{itemize}

This hypothesis is satisfied for many scattering maps, e.g.~specular reflection\footnote{Here $\theta(w)=\pi-2\arcsin(w)$ and thus condition (A) holds.} or the scattering in the muffin-tin Coulomb potential $V(\vecq)=\alpha \max(\|\vecq\|^{-1}-1,0)$ with $\alpha\notin \{-2E,0\}$, where $E$ denotes the total energy, cf.~\cite{partII}.

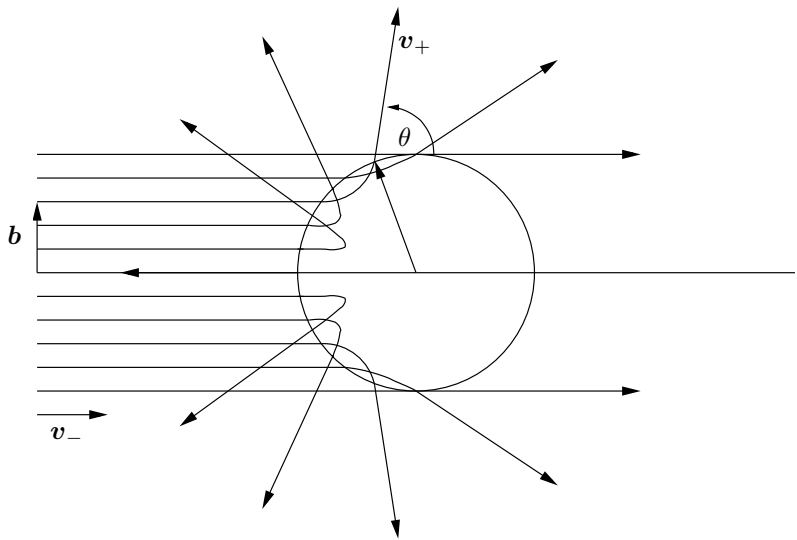
\begin{figure}
\begin{center}
\begin{minipage}{0.9\textwidth}
\unitlength0.1\textwidth
\begin{picture}(10,6.3)(0,0)
\put(0.3,0){\includegraphics[width=0.9\textwidth]{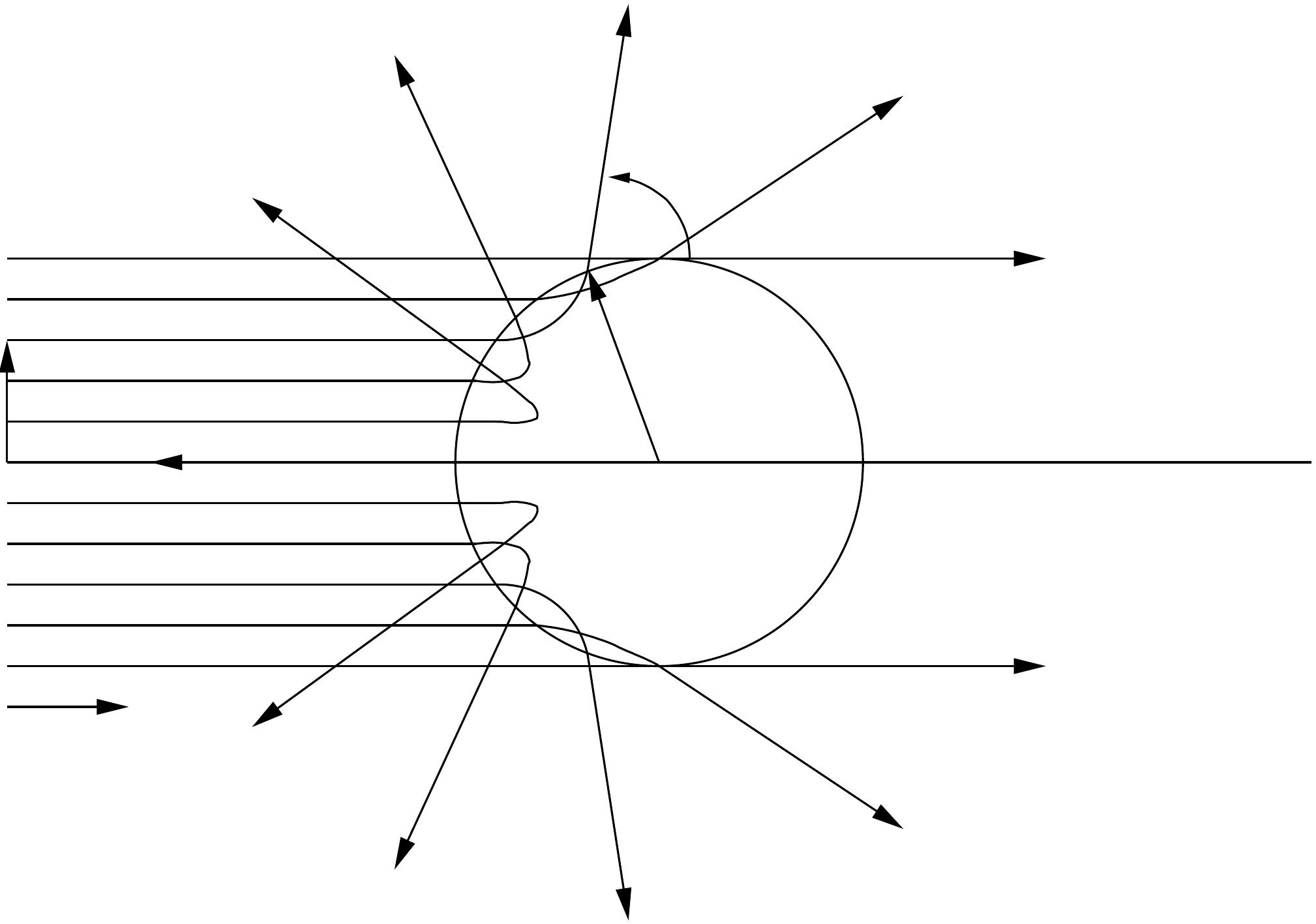}}
\put(0,3.5){$\vecb$} 
\put(0.5,1.2){$\vecv_-$} 
\put(4.6,4.65){$\theta$}
\put(4.6,5.8){$\vecv_+$} 
\end{picture}
\end{minipage}
\end{center}
\caption{Illustration of a scattering map satisfying Hypothesis (A).} \label{fig4}
\end{figure}

An inductive argument shows that there is a sequence $(\vecw_n)_{n\in\NN}$ in $\UB$, so that the impact parameter $\vecb_n$,  exit velocity $\vecv_n$ and exit parameter $\vecs_n$ at the $n$th collision are given by the frame-independent formulas
\begin{equation}\label{roti}
\vecv_n= \vece_1 R_n^{-1} , \qquad \vecb_n= (0,\vecw_n)R_{n-1}^{-1}
\qquad \vecs_n= (0,\vecw_n) R_n^{-1}
\end{equation}
where
\begin{equation}\label{Rn}
R_n := R(\vecv_0) S(\vecw_1)\cdots S(\vecw_n) .
\end{equation}
Here $R:\S_1^{d-1}\to\SO(d)$ is smooth  up to finitely many singular points, such that
$\vecv R(\vecv)=\vece_1$ for all $\vecv\in\S_1^{d-1}$. For an example see footnote 3 on p.~1968 of \cite{partI}.

We can now express position and velocity at time $t>0$ as\footnote{The $O(r \nu(t))$-error is simply due to the fact that we have not included the jumps of position at each scattering. In the case of specular reflections, all formulas are exact.}
\begin{equation}\label{qdef}
\vecq(t) =  \vecq_{\nu(t)} + (t-\tau_{\nu(t)})\vecv(t) +O(r \nu(t)),\qquad \vecv(t)=\vecv_{\nu(t)} , 
\end{equation}
where 
\begin{equation}
\tau_n := \sum_{j=1}^n t_j,\quad \tau_0:=0, 
\end{equation}
is the time to the $n$th collision, $t_j$ is the $j$th intercollision time,
\begin{equation}
\nu(t) := \max\{ n\in\ZZ_{\geq 0} : \tau_n \leq t \} 
\end{equation}
is the number of collisions within time $t$,
\begin{equation}
\vecq_n := \sum_{j=1}^n t_j \vecv_{j-1}  
\end{equation}
is the particle location at the $n$th collision\footnote{Again, this is up to an error of order $O(r n)$.} and 
\begin{equation}
\vecv_n = R(\vecv_0) S(\vecw_1) \cdots S(\vecw_n) \vece_1
\end{equation}
is the velocity after the $n$th collision as calculated in \eqref{roti}.

In the macroscopic coordinates \eqref{macrocoo}, the above translates to
\begin{equation}
\vecQ(t) =  \vecQ_{\scrV(t)}  + (t-\scrT_{\scrV(t)} )\vecV(t) +O(r^d \scrV(t)),\qquad \vecV(t)=\vecV_{\scrV(t)} 
\end{equation}
where $\vecQ_n=r^{d-1}\vecq_n$, $\vecV_n=\vecv_n$, $\scrT_n=r^{d-1}\tau_n$ and 
\begin{equation}\label{scrVdef}
\scrV(t) := \nu(r^{-(d-1)} t)=\max\{ n\in\ZZ_{\geq 0} : \scrT_n \leq t \} .
\end{equation}

\section{A refined Stosszahlansatz}\label{sec:beyond}

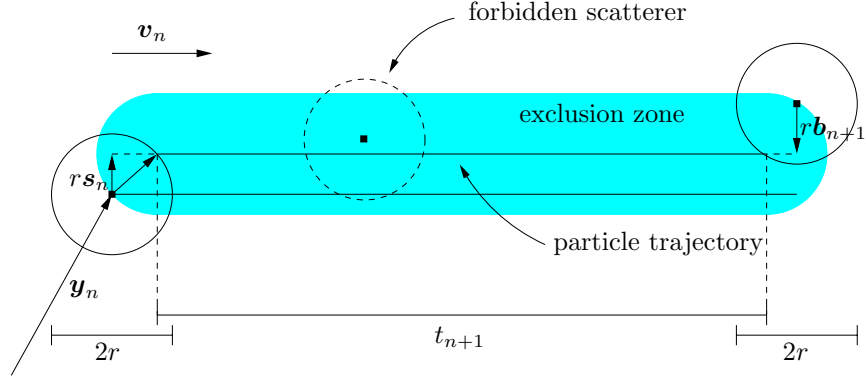
\begin{figure}
\begin{center}
\begin{minipage}{0.9\textwidth}
\unitlength0.1\textwidth
\begin{picture}(10,4.5)(0,0)
\put(0,0){\includegraphics[width=\textwidth]{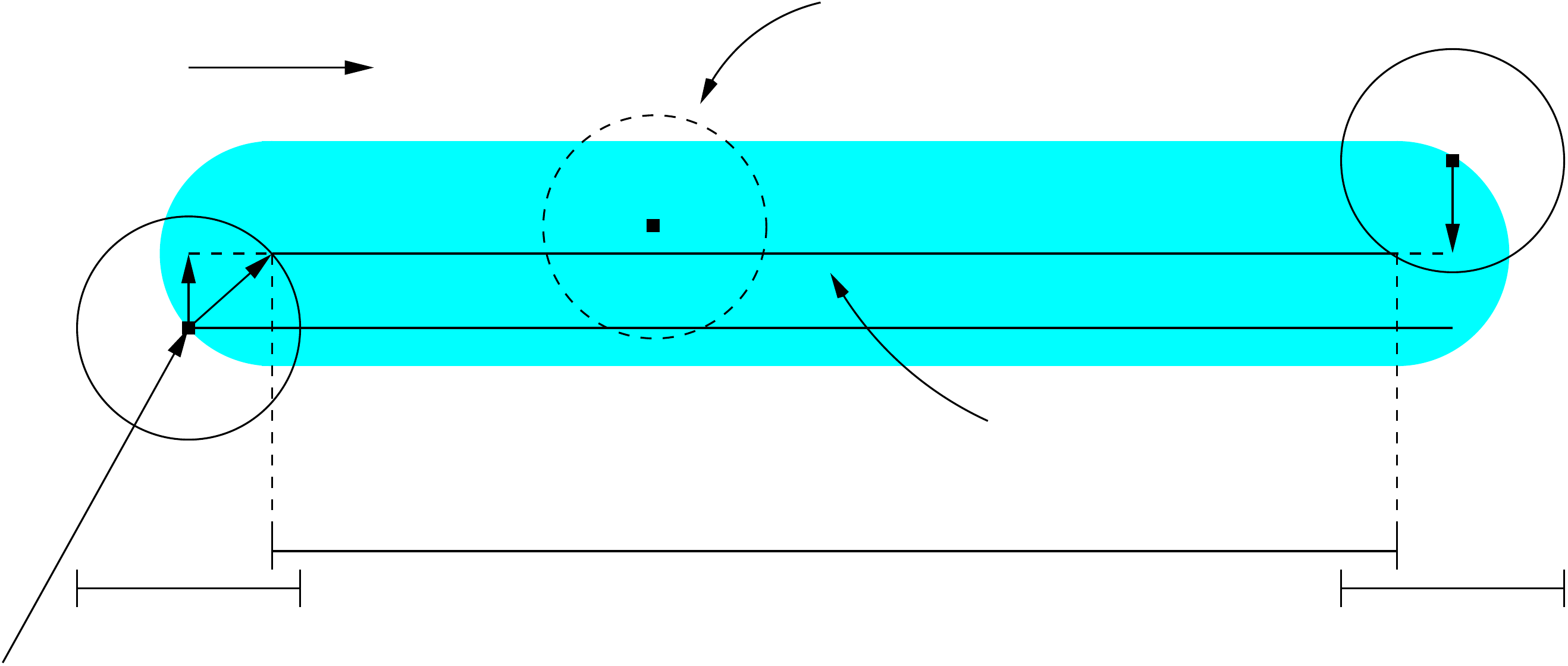}}
\put(0.7,2.25){$r\vecs_n$} 
\put(0.7,1){$\vecy_n$}
\put(1.5,4){$\vecv_n$} 
\put(1,0.2){$2r$}
\put(9.1,0.2){$2r$}
\put(9.3,2.85){$r\vecb_{n+1}$}
\put(5,0.4){$t_{n+1}$}
\put(5.4,4.2){forbidden scatterer}
\put(6.4,1.5){particle trajectory}
\put(6,3){exclusion zone}
\end{picture}
\end{minipage}
\end{center}
\caption{Intercollision flight in the Lorentz gas between the $n$th and $(n+1)$st collision. The exclusion zone is a cylinder of radius $r$ with spherical caps.} \label{fig5}
\end{figure}

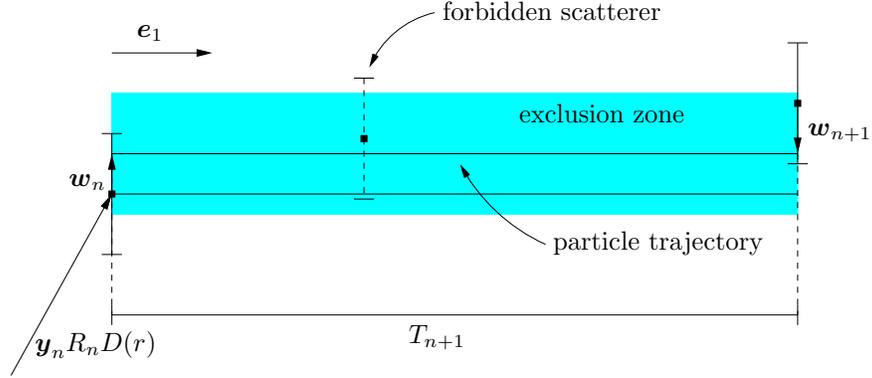
\begin{figure}
\begin{center}
\begin{minipage}{0.9\textwidth}
\unitlength0.1\textwidth
\begin{picture}(10,4.5)(0,0)
\put(0,0){\includegraphics[width=\textwidth]{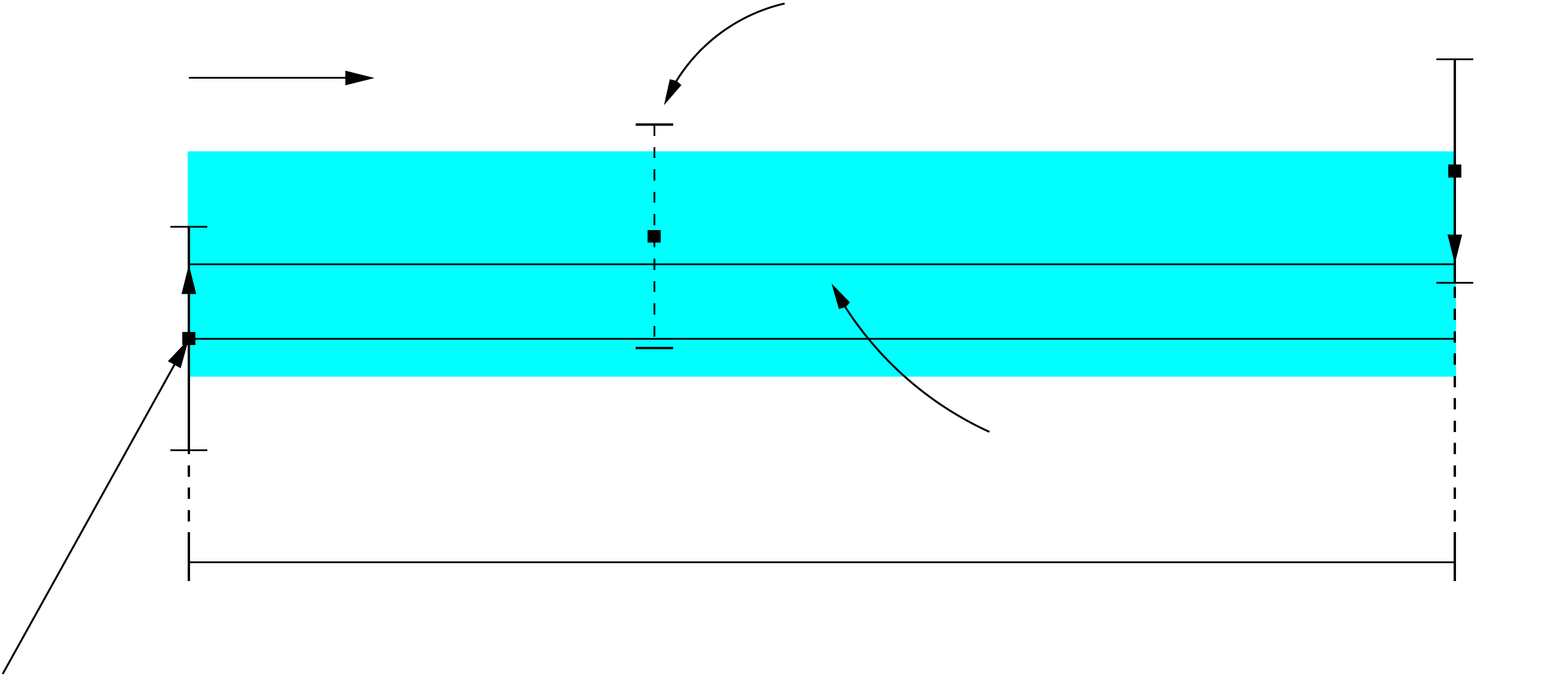}}
\put(0.7,2.25){$\vecw_n$} 
\put(0.3,0.3){$\vecy_n R_n D(r)$}
\put(1.5,4){$\vece_1$} 
\put(9.4,2.85){$\vecw_{n+1}$}
\put(4.7,0.4){$T_{n+1}$}
\put(5.1,4.2){forbidden scatterer}
\put(6.4,1.5){particle trajectory}
\put(6,3){exclusion zone}
\end{picture}
\end{minipage}
\end{center}
\caption{The intercollision flight in Fig.~\ref{fig5} after applying the linear map $R_n D(r)$ with $r$ very small. The exclusion zone is now approximately a cylinder with flat caps.} \label{fig6}
\end{figure}

We will now investigate the particle trajectory corresponding to random initial conditions $(\vecQ_0,\vecV_0)$ and outline a strategy to establish the convergence to a random flight process in the Boltzmann-Grad limit.\footnote{We assume here that $(\vecQ_0,\vecV_0)$ is distributed according to a fixed, absolutely continuous probability measure $\Lambda$ on $\T^1(\RR^d)$. One can, of course also prepare the initial particle cloud on smaller scales. For example take $(\vecq_0,\vecv_0)=(r^{-(d-1)}\vecQ_0,\vecV_0)$ random with respect to a fixed absolutely continuous $\Lambda$. In the case of the periodic and the quasicrystal Lorentz gas \cite{partI,partII,quasi} we are even able to consider more singular $\Lambda$: Fix $\vecq_0$ and only take $\vecv_0$ random according to an absolutely continuous measure on the unit sphere. In this case, we have convergence for {\em every} $\vecq_0$, with the same limit distribution for {\em almost every} $\vecq_0$.}  We will, for now, keep the scatterer configuration $\scrP$ general, and discuss in later sections examples of $\scrP$ which allow a rigorous treatment.

Let us focus on the $n$th and $(n+1)$st collision and consider a parallel beam of particles with given velocity $\vecv_{n-1}$ that hit a scatterer located at $\vecy_n$ with a certain intensity distribution $\lambda$ in the impact parameter $\vecw_n$ (Figure \ref{fig5}).\footnote{The measure $\lambda$ will of course depend on the history of the particle beam, and in particular on $r$, but let us assume here for the sake of argument that $\lambda$ is a {\em fixed} Borel probability measure on $\UB$. A key part of the paper \cite{partII} deals with the problem of $r$-dependent measures in the setting of the periodic Lorentz gas, by obtaining uniform estimates over families of $\lambda$.} The task is now to calculate the probability of hitting the next scatterer in a small time interval around $t_{n+1}$ with impact parameter near $\vecw_{n+1}$. Recall that we expect $t_{n+1}$ to be of order $r^{-(d-1)}$, and it is natural to set $T_n=r^{d-1}t_n$. 
We now first shift our coordinate system by $-\vecy_n-r(\vecs_n+\vecv_n\sqrt{1-\|\vecs_n\|^2})$ so the left center of the cylinder is now at the origin, then rotate our coordinate system by $R_n\in\SO(d)$, so that the outgoing velocity $\vecv_n$ becomes $\vece_1$, cf.~\eqref{roti}, and finally apply the linear transformation given by the matrix
\begin{equation}
D(r) =
\begin{pmatrix}
r^{d-1} & \vecnull \\ \trans\vecnull  & r^{-1} 1_{d-1} 
\end{pmatrix}
\end{equation}
which rescales the length units along and perpendicular to the cylinder.
Note that the caps of the cylinder become flat as $r\to 0$, cf.~Fig.~\ref{fig6}. In particular
\begin{equation}
\begin{split}
r\big(\vecs_n+\vecv_n\sqrt{1-\|\vecs_n\|^2}\big) R_n D(r) 
& = r\big((0,\vecw_n)+(1,\vecnull)\sqrt{1-\|\vecw_n\|^2}\big) D(r) \\
& = (0,\vecw_n)+ O(r^d).
\end{split}
\end{equation}
The rotation matrix $R_n$ is, by \eqref{Rn}, given by $R_n=R_{n-1} S(\vecw_n)$ where $R_{n-1}$ is fixed (since $\vecv_{n-1}$ is assumed fixed in this discussion). For $\vecw_n$ random according to $\lambda$, we are interested in the probability that the particle hits the next scatterer at a time $T_{n+1}$ in the interval $A=\,]\xi,\xi+d\xi[$ and with impact parameter $\vecw_{n+1}$ in some box $B\subset\UB$. This probability is, for small $r$, approximately\footnote{This approximation is justified, if the limit distribution is continuous in $\xi$.\label{foot1}} equal to the probability that the random point set
\begin{equation}\label{theta1}
\widetilde\Theta_r(\vecy_n) = (\scrP-\vecy_n) R_{n-1} S(\vecw_n) D(r) - (0,\vecw_n)
\end{equation}
does not intersect the cylinder 
$\fZ(\xi)=]0,\xi[\times\UB$
and has (at least\footnote{We assume that, in the limit $r\to 0$, the probability of having one point in a small set is approximately the same as the probability of having one ore more points. As in footnote \ref{foot1}, this is justified, if the limit distribution is continuous in $\xi$.\label{foot2}}) one point in the box
$A\times B$.
Our general objective is therefore to try to prove that there is a random point process\footnote{Throughout this paper, we will represent random point processes as random point sets.} $\Theta(\vecy)$ in $\RR^d$ and a random variable $\vech\in\UB$ distributed according to $\lambda$ such that, for every fixed $\vecy\in\scrP$, 
\begin{equation}\label{limitTheta}
\widetilde\Theta_r(\vecy) \xrightarrow[r\to 0]{} \widetilde\Theta(\vecy):=\Theta(\vecy) - (0,\vech)
\end{equation}
in finite-dimensional distribution. This means that for any $k\in\NN$, $\scrA_1,\ldots,\scrA_k\subset\RR^d$ bounded with boundary of measure zero and $n_1,\ldots,n_k\in\ZZ_{\geq 0}$, we have
\begin{equation}\label{statem00}
\lim_{r\to 0} \PP\big(\#(\widetilde\Theta_r(\vecy)\cap\scrA_i) = n_i \,\forall i)  = \PP\big( \#(\widetilde\Theta(\vecy)\cap \scrA_i)=n_i \, \forall i \big) .
\end{equation}

It is crucial that $\Theta(\vecy)$ and $\vech$ are independent, and that $\Theta(\vecy)$ is independent of the choice of $\lambda$ and $R_n$. We conclude that, if the convergence in \eqref{limitTheta} indeed holds in finite-dimensional distribution (as we are dealing with only two test sets, $\fZ(\xi)$ and $A\times B$, convergence in two-dimensional distribution is in fact sufficient)  then the probability that the particle hits the next scatterer at a time $T_{n+1}\in A$ and with impact parameter $\vecw_{n+1}\in B$, is in the limit $r\to 0$ given by
\begin{equation}\label{sat8}
\PP\big( \widetilde\Theta(\vecy_n)\cap\fZ(\xi)=\emptyset ,\;  \#(\widetilde\Theta(\vecy_n)\cap (A\times B)) =1 \big) .
\end{equation}
In some instances, $\Theta(\vecy)$ will not depend on the scatterer location $\vecy$, for example when $\scrP$ is a realisation of a Poisson process or a Euclidean lattice, as we shall see below. If $\Theta(\vecy)$ does depend on the scatterer location, the hope is that this dependence is ``mild,'' in the sense that there exists a probability space $(\Sigma,\scrF,\mm)$ and a map
\begin{equation}\label{iotamap}
\iota: \scrP \to \Sigma,\qquad \vecy \mapsto \iota(\vecy),
\end{equation}
so that $\Theta(\vecy)$ depends only on the value of $\iota(\vecy)$. We will call $\iota(\vecy)$ the {\em colour} of $\vecy$, and consider the colourised scatterer configuration, 
\begin{equation}\label{colours}
\{ (\vecy, \iota(\vecy)) : \vecy\in\scrP \} \subset \RR^d\times\Sigma.
\end{equation}
We assume furthermore that the colour in \eqref{colours} is distributed according to the probability measure $\mm$ on $\Sigma$, in the sense that (cf.\ \eqref{density})
for any bounded $\scrD\subset\RR^d$ with $\vol_{\RR^d}(\scrD)>0$, $\vol_{\RR^d}(\partial\scrD)=0$ and any measurable set $B\subset\Sigma$ with $\mm(\partial B)=0$,
\begin{equation}\label{density2}
\lim_{R\to\infty} \frac{\#\{\vecy\in\scrP\cap R\scrD : \iota(\vecy)\in B\}}{\vol_{\RR^d}(R\scrD)} =\mm(B).
\end{equation}
Let us define $\Omega:=\Sigma\times\UB$ as the product space of colour and impact parameters, with probability measure $\pp=\mm\times v_{d-1}^{-1} \vol_{\RR^{d-1}}$.
Instead of \eqref{limitTheta}, we must now consider the convergence for the corresponding colourised point processes. Once we understand the colourised limit, we can compute the limit distribution for the probability of emerging from a scatterer with a given colour and exit parameter $\omega_n$, and hitting the next scatterer at time $T_n\in\;]\xi,\xi+d\xi[$ with colour and impact parameter $\omega_{n+1}\in B \subset\Omega$. We denote this probability by
\begin{equation}\label{sat9}
\int_B  k(\omega_n,\xi,\omega) \,d\xi \,\pp(d\omega) ,
\end{equation}
which defines the {\em transition kernel} $k(\omega',\xi,\omega)$.
The conclusion of the above heuristics is now that the particle trajectory
\begin{equation}\label{traj}
\Xi_r: t \mapsto (\vecQ(t),\vecV(t)),
\end{equation}
with random initial condition $(\vecQ_0,\vecV_0)$ distributed according to some absolutely continuous measure $\Lambda$ on $\T^1(\RR^d)$, converges in the Boltzmann-Grad limit to the continuous-time random flight process $\Xi(t)$ in $\T^1(\RR^d)$ defined as follows. 

Consider the sequences of random variables $\uxi=(\xi_n)_{n\in\NN}$ and $\ueta=(\eta_n)_{n\in\NN}$ defined by the Markov chain
\begin{equation}\label{Mchaindef}
n \mapsto (\xi_n, \eta_n) 
\end{equation}
with state space $\RR_{>0}\times\Omega$ and transition probability $(n\geq 2)$
\begin{equation}\label{trapro}
\PP((\xi_n, \eta_n)\in A \mid \xi_{n-1}, \eta_{n-1}) = \int_A  k(\eta_{n-1},\xi,\omega) \,d\xi \,\pp(d\omega) ,
\end{equation}
where the transition kernel $k(\omega',\xi,\omega)$ is defined by \eqref{sat9}. The initial distribution is
\begin{equation}
\PP((\xi_1, \eta_1)\in A ) = \int_A  K(\xi,\omega) \,d\xi \,\pp(d\omega) ,
\end{equation}
where
\begin{equation}\label{bigK}
K(\xi,\omega) := \frac{1}{\,\xibar\,}\int_\xi^\infty \int_\Omega  k(\omega',\xi',\omega) \,\pp(d\omega')\,d\xi' .
\end{equation}
The time-reversibility of the underlying microscopic dynamics (for every fixed $r>0$) implies that the transition kernel $k$ is symmetric, i.e.
\begin{equation}\label{ksymm}
k(\omega,\xi,\omega') = k(\omega',\xi,\omega).
\end{equation}
Because the transition probability \eqref{trapro} is independent of $\xi_{n-1}$, the chain $n \mapsto \eta_n$ is also Markovian, with transition probability
\begin{equation}\label{trapro2}
\PP(\eta_n \in A \mid \eta_{n-1}) = \int_A  \int_0^\infty k(\eta_{n-1},\xi,\omega) \,d\xi \,\pp(d\omega) .
\end{equation}
The stationary measure for this Markov chain is $\pp$, and the distribution of free path lengths with respect to this measure is defined as
\begin{equation}\label{FPL}
\Psi_0(\xi):= \int_\Omega\int_\Omega k(\omega',\xi,\omega) \,\pp(d\omega)\,\pp(d\omega').
\end{equation}

Let us write $\eta_n=(\chi_n,\vech_n)$, where $\chi_n\in\Sigma$ is the colour and $\vech_n\in\UB$ the impact parameter. 
In analogy with the deterministic setting \eqref{qdef}--\eqref{scrVdef}, we define the random variables
\begin{equation}\label{Tdef22}
\scrT_n^\BG := \sum_{j=1}^n \xi_j,\quad \scrT_0^\BG:=0, 
\end{equation}
\begin{equation}
\scrV^\BG(t) := \max\{ n\in\ZZ_{\geq 0} : \scrT_n^\BG \leq t \} ,
\end{equation}
\begin{equation}
\vecQ_n^\BG := \vecQ_0+\sum_{j=1}^n \xi_j \vecV_{j-1}^\BG
,\qquad 
\vecV_n^\BG := R(\vecV_0) S(\vech_1) \cdots S(\vech_n) \vece_1 ,
\end{equation}
\begin{equation}\label{Vdef22}
\vecQ^\BG(t) :=  \vecQ_{\scrV^\BG(t)}^\BG + (t-\scrT_{\scrV^\BG(t)}^\BG)\vecV^\BG(t),
\qquad
\vecV^\BG(t):=\vecV_{\scrV^\BG(t)}^\BG .
\end{equation}
The notation ``BG'' stands for {\em Boltzmann-Grad limit} and is used to differentiate from the deterministic counterparts \eqref{qdef}--\eqref{scrVdef}.
Note that none of the above depend explicitly on colour. The hidden variable ``colour'' is needed to make \eqref{Mchaindef} a Markov chain.
The random flight process $\Xi$ is thus defined as
\begin{equation}\label{Xidef2}
t\mapsto \Xi(t) := \big(\vecQ^\BG(t),\vecV^\BG(t)\big) .
\end{equation}
The convergence of the random process $\Xi_r$ in \eqref{traj} to $\Xi$ answers in particular our question \eqref{Lorentz-limit}, since the former implies the convergence in \eqref{Lorentz-limit} with $L^t$ defined by
\begin{equation}\label{Ltdef}
\int_\scrA L^t f(\vecQ,\vecV) \, d\vecQ\,d\vecV = \PP( \Xi(t)\in\scrA) .
\end{equation}
Here $f=\Lambda'$ is the Radon--Nikodym derivative of $\Lambda$.

\section{A generalised Boltzmann equation}\label{sec:generalized}

This limiting process $\Xi(t)$ defined in \eqref{Xidef2} is in general not a continuous-time Markov process,\footnote{A consequence of this fact is that the family of operators $L^t$ in \eqref{Ltdef} does not form a semigroup, i.e., $L^t L^s= L^{t+s}$ does not hold for all $s,t>0$.} but can be turned into one by extending the state space as follows. We define the time until the next scattering by
\begin{equation}
T^\BG(t):= \scrT_{\scrV^\BG(t)+1}^\BG-t ,
\end{equation}
the colour of the next scatterer by
\begin{equation}
\chi^\BG(t):= \chi_{\scrV^\BG(t)+1} ,
\end{equation}
and the exit velocity of the next scattering by
\begin{equation}
\vecV_+^\BG(t):=\vecV_{\scrV^\BG(t)+1}^\BG .
\end{equation}
The process
\begin{equation}
t\mapsto \widetilde\Xi(t):=\big(\vecQ^\BG(t),\vecV^\BG(t),T^\BG(t),\chi^\BG(t),\vecV_+^\BG(t)\big) 
\end{equation}
is now a Markov process with state space $\T^1(\RR^d)\times\RR_{>0}\times\Sigma\times\US$ and backward equation\footnote{This equation is also known as {\em Fokker-Planck-Kolmogorov equation}.}
\begin{equation}\label{FPK}
\begin{cases}
(\partial_t +\vecV\cdot\partial_\vecQ-\partial_\xi) f_t(\vecQ,\vecV,\xi,\chi,\vecV_+) 
= [\scrC f_t](\vecQ,\vecV,\xi,\chi,\vecV_+) & \\[10pt]
\displaystyle
\lim_{t\to 0} f_t(\vecQ,\vecV,\xi,\chi,\vecV_+) = \Lambda'(\vecQ,\vecV)\, K(\xi,\omega)\, \sigma(\vecV,\vecV_+),&
\end{cases}
\end{equation}
with $K(\xi,\omega)$ as in \eqref{bigK} and  the collision operator $\scrC$ is defined by
\begin{multline}
 [\scrC f](\vecQ,\vecV,\xi,\chi,\vecV_+) \\
 = \sigma(\vecV,\vecV_+) \int_\US \int_\Sigma f(\vecQ,\vecV',0,\chi',\vecV) \,k(\omega',\xi,\omega) \,d\mm(\chi') \, d\vecV' ,
\end{multline}
where
\begin{equation}
\omega':=(\chi',\vecs(\vecV',\vecV)R(\vecV)),\qquad \omega:=(\chi,\vecb(\vecV,\vecV_+)R(\vecV)) .
\end{equation}

A stationary solution of eq.~\eqref{FPK} is given by
\begin{equation}\label{FPKstat}
f_t(\vecQ,\vecV,\xi,\chi,\vecV_+) = K(\xi,\omega)\, \sigma(\vecV,\vecV_+),
\end{equation}
which corresponds to $\Lambda=$ Liouville measure. To see this, note that the left hand side of the first line in \eqref{FPK} is
\begin{equation} \label{lhs}
\sigma(\vecV,\vecV_+) \, \xibar^{-1} \int_\Omega k(\omega',\xi,\omega)\, d\pp(\omega')  .
\end{equation}
Furthermore, we have
\begin{equation*}
\xibar K(0,\omega') = \int_0^\infty \int_\Omega k(\omega'',\xi,\omega')\, d\xi\, d\pp(\omega'') 
= \int_0^\infty \int_\Omega k(\omega',\xi,\omega'')\, d\xi\, d\pp(\omega'') 
= 1.
\end{equation*}
The right hand side of the first line in \eqref{FPK} therefore equals, in view of \eqref{dcs},
\begin{multline}\label{rhs}
\sigma(\vecV,\vecV_+) \int_\US \int_\Sigma \sigma(\vecV',\vecV)\, K(0,\omega')\,k(\omega',\xi,\omega) \,d\mm(\chi') \, d\vecV' 
\\
= \sigma(\vecV,\vecV_+) \, \xibar^{-1} \int_\US \int_\Sigma k(\omega',\xi,\omega) \,d\pp(\omega') ,
\end{multline}
which equals \eqref{lhs}
This shows that \eqref{FPKstat} is indeed a stationary solution of \eqref{FPK}.

Let us now illustrate the above programme with a number of examples, where all or part of the heuristics can be made rigorous. The principal questions we would like to answer, for a given scatterer configuration $\scrP$, are: {\em Does the limit \eqref{Lorentz-limit} exist?} {\em What is the limit process $\Theta(\vecy)$?} {\em What is the transition kernel $k(\omega',\xi,\omega)$?} 

We begin with the classic setting where $\scrP$ is a typical realisation of a Poisson process and will show how the generalised linear Boltzmann equation \eqref{FPK} reduces to the original.

\section{Random scatterer configuration}\label{sec:Poisson}

The Poisson process $\Theta=\Theta_{\text{\rm Poisson}}$ in $\RR^d$ with intensity $\nbar=1$ is characterised by the property that for any collection of bounded, pairwise disjoint Borel sets $\scrA_1,\ldots,\scrA_k$ and integers $n_1,\ldots,n_k\geq 0$, 
\begin{equation}
\PP(\# (\Theta\cap\scrA_i)=n_i \;\forall i ) 
= \prod_{i=1}^k \frac{(\vol_{\RR^d}(\scrA_i))^{n_i}}{n_i!} \e^{-\vol_{\RR^d}(\scrA_i)} .
\end{equation}
We will assume in this section that $\scrP$ is a fixed realisation of a Poisson process. In a seminal paper, Boldrighini, Bunimovich and Sinai \cite{Boldrighini83} have shown that the limit \eqref{Lorentz-limit} exists almost surely and is given by the linear Boltzmann equation \eqref{LBeq}. 

\begin{thm}[Boldrighini, Bunimovich and Sinai, 1983 \cite{Boldrighini83}]\label{BBS}
The convergence in \eqref{Lorentz-limit} holds for a typical realisation $\scrP$ of a Poisson process, and $f_t=L^t f$ satisfies the linear Boltzmann equation \eqref{LBeq}.
\end{thm}

This result was previously proved by Gallavotti \cite{Gallavotti69} on average for random $\scrP=\Theta_{\text{\rm Poisson}}$, and by Spohn \cite{Spohn78} for more general random scatterer configurations and scattering potentials.

In the present setting, the limit process $\Xi(t)$ is in fact already a continuous time Markov process and the extension to $\widetilde\Xi(t)$ is not necessary. Nevertheless it is instructive to see how the backward equation \eqref{FPK} reduces to the linear Boltzmann equation \eqref{LBeq}. 

A review of the arguments used in \cite{Boldrighini83} shows that the convergence \eqref{limitTheta} holds in finite-dimensional distribution for almost all $\scrP$ with limit $\Theta(\vecy)=\Theta_{\text{\rm Poisson}}$ and thus, by the translation invariance of the Poisson process, $\widetilde\Theta(\vecy)=\Theta_{\text{\rm Poisson}}$. The limiting point process is evidently independent of $\vecy$, and we may paint all scatterers in the same colour. That is, $\Sigma$ is the space of one element. We can thus identify $\Omega$ with $\UB$ and set $\pp(d\vecw)=v_{d-1}^{-1} d\vecw$. The Poisson distribution yields in \eqref{sat8} the transition kernel 
\begin{equation}
k(\omega',\xi,\omega)=\xibar^{-1}\, \e^{-\xi/\xibar} ,\qquad
K(\xi,\omega)=\xibar^{-1}\, \e^{-\xi/\xibar}  .
\end{equation}
The ansatz
\begin{equation}
f_t(\vecQ,\vecV,\xi,\chi,\vecV_+) =g_t(\vecQ,\vecV)\,  \sigma(\vecV,\vecV_+) \, \xibar^{-1}\, \e^{-\xi/\xibar} 
\end{equation}
in the backward equation \eqref{FPK} of $\widetilde\Xi(t)$ shows that, after a separation of variables, the function $g_t(\vecQ,\vecV)$ is a solution of the linear Boltzmann equation \eqref{LBeq}. More directly, one can show that $\Xi(t)$ is Markov, and that  the linear Boltzmann equation is the backward equation of $\Xi(t)$.

\section{Periodic scatterer configuration}\label{sec:periodic}

The opposite extreme of a random scatterer configuration is a perfectly periodic point set $\scrP$. We assume in this section that $\scrP$ is a Euclidean lattice $\scrL$ of covolume one. More general periodic scatterer configurations are considered as a special case in the framework of quasicrystals, cf.~Section \ref{sec:quasicrystals}.

\begin{thm}[Marklof and Str\"ombergsson, 2008 \cite{partII}]\label{thm3}
The convergence in \eqref{Lorentz-limit} holds for every Euclidean lattice $\scrP=\scrL$ of covolume one, where $L^t$ is independent of the choice of $\scrL$.
\end{thm}

The main result of \cite{partII} is in fact more general: It extends to the convergence in distribution of the random process $\Xi_r$ in \eqref{traj} to $\Xi$. The proof of Theorem \ref{thm3} turns the heuristics of Section \ref{sec:beyond} into a rigorous argument. Let us describe some of the key objects.

Every Euclidean lattice of covolume one can be written as $\scrL =\ZZ^d M$
for some $M\in\SLdR$. Since the stabiliser of $\ZZ^d$ under right multiplication by $G=\SLdR$ is the subgroup $\Gamma=\SLdZ$, one can show that there is a bijection
\begin{equation}\label{SOL}
\begin{split}
\GamG & \xrightarrow{\sim} \{ \text{Euclidean lattices of covolume one} \} \\
\Gamma M & \mapsto \ZZ^d M .
\end{split}
\end{equation}
It is a well known fact that any fundamental domain of $\Gamma=\SLdZ$ has finite Haar measure in $G=\SLdR$. This implies that there is a unique probability measure $\mu$ on $\GamG$ invariant under the natural $G$-action (which is multiplication from the right). We define a random point process in $\RR^d$ by setting $\Theta_{\text{\rm lattice}} = \ZZ^d M$
with $M$ random in $\GamG$ according to $\mu$ and the above identification \eqref{SOL} of $\GamG$ and the space of lattices. We will call $\Theta_{\text{\rm lattice}}$ a {\em random lattice}.

The following theorem says that, for any fixed $\scrP=\scrL$ the convergence in \eqref{limitTheta} holds with $\Theta=\Theta_{\text{\rm lattice}}$. Note that by translational invariance of $\scrL$, all point processes in \eqref{limitTheta} are independent of $\vecy$, and we will write in the following $\widetilde\Theta_r$ instead of $\widetilde\Theta_r(\vecy)$.

\begin{thm}[\cite{partI}]\label{thm13.1}
Let $\lambda$ be an absolutely continuous probability measure on $\UB$, $\scrA_1,\ldots,\scrA_k\subset\RR^d$ bounded with boundary of measure zero and $n_1,\ldots,n_k\in\ZZ_{\geq 0}$. Then
\begin{equation}\label{statem}
\lim_{r\to 0} \PP\big(\#(\widetilde\Theta_r\cap\scrA_i) = n_i \,\forall i\big)  = \PP\big( \#((\Theta_{\text{\rm lattice}}-(0,\vech))\cap \scrA_i)=n_i \, \forall i \big) .
\end{equation}
\end{thm}

This theorem is a consequence of equidistribution of large spheres on $\GamG$: 

\begin{thm}[\cite{partI}]\label{thm13.2}
For any $M\in\GamG$, any bounded continuous $f: \UB\times \GamG \to \RR$ and any absolutely continuous probability measure $\lambda$ on $\UB$, 
 \begin{equation}\label{eq13.2i}
\lim_{r\to 0}\int_{\UB} f(\vecw, M S(\vecw) D(r))\, d\lambda(\vecw)  =
\int_{\UB} \int_{\GamG} f(\vecw,M) \,d\mu(M)\, d\lambda(\vecw).
\end{equation}
\end{thm}

Theorem \ref{thm13.1} is derived from Theorem \ref{thm13.2} by choosing in \eqref{eq13.2i} as test function $f$ the characteristic function of the set
\begin{equation}\label{0set}
\big\{ (\vecw,M) \in \UB \times \GamG :  \#\big((\ZZ^d M - (0,\vecw )\cap\scrA_i \big)= n_i \,\forall i\big\} .
\end{equation}
This choice does of course not produce a continuous $f$, but one can show that \eqref{0set} has boundary of measure zero in $\UB\times\GamG$, and thus the characteristic function can be approximated sufficiently well by continuous functions. Details of this technical argument can be found in \cite{partI}, Sections 5 and 6.

Since the limit process $\Theta_{\text{\rm lattice}}$ is independent of $\vecy$ there is no need for colour (as in the Poisson setting), and we again identify $\Omega$ with $\UB$, and set $\pp(d\vecw)=v_{d-1}^{-1} d\vecw$. 
In order to work out the transition kernel $k(\vecw',\xi,\vecw)$ in \eqref{sat9}, set $X=\GamG$ and define the subspace
\begin{equation}
X(\vecy)=\{ M\in X : \vecy\in\ZZ^d M \}
\end{equation}
of those lattices (of covolume one) that contain a given $\vecy\in\RR^d$.
In \cite{partI} we construct a probability measure $\nu_\vecy$ on $X(\vecy)$ so that 
\begin{equation}
d\mu(M) = d\nu_\vecy(M)\, d\vecy.
\end{equation}
With this, we can infer that
\begin{equation}\label{1st}
k(\vecw',\xi,\vecw) = \xibar^{-1} \nu_\vecy\big(\big\{  M\in X(\vecy) : \ZZ^d M \cap (\fZ(\xi)+(0,\vecw'))=\emptyset \big\}\big)
\end{equation}
where $\vecy=(\xi,\vecw'-\vecw)$. For an explicit description of the $\nu_\vecy$-measure of the above set, see \cite{partIV}, Section 2.2.
In dimension $d=2$, when $\scrB_1^1=\,]\!-1,1[\,$, eq.~\eqref{1st} can be used to calculate an explicit formula for the transition kernel. We have \cite{partIII}
\begin{equation}\label{twodidi}
k(\vecw',\xi,\vecw)=  \frac{12}{\pi^2}\, \Upsilon\bigg( 1+ \frac{\xi^{-1} - \max(|\vecw|,|\vecw'|)-1}{|\vecw-\vecw'|}\bigg)
\end{equation}
with
\begin{equation}
\Upsilon(x)=
\begin{cases} 
0 & \text{if }x\leq 0\\
x & \text{if }0<x<1\\
1 & \text{if }x\geq 1.
\end{cases}
\end{equation}
For independent derivations of Formula \eqref{twodidi} that do not employ eq.~\eqref{1st} but a more direct approach based on Farey dissections, see  Bykovskii and Ustinov \cite{Bykovskii09} and Caglioti and Golse \cite{Caglioti10}.

There are no such formulas in higher dimension, although \eqref{1st} can be used to extract information to obtain asymptotics for $\xi\to 0$ and $\xi\to\infty$, cf.~\cite{partIV}. We have in particular 
\begin{equation}\label{Kbound}
\frac{1-2^{d-1} \xibar^{-1}\xi}{\zeta(d)\xibar} \leq k(\vecw',\xi,\vecw) \leq \frac{1}{\zeta(d)\xibar},
\end{equation}
and so for small $\xi$ this implies $k(\vecw',\xi,\vecw) = (\zeta(d)\xibar)^{-1}+O(\xi)$.
Here $\zeta(d)$ is the Riemann zeta function and $\zeta(d)^{-1}$ is the relative density of primitive lattice points in $\ZZ^d$.
Compare \eqref{Kbound} with the result for the Poisson process (Section \ref{sec:Poisson}):
\begin{equation}
k_\text{Poisson}(\vecw',\xi,\vecw) = \xibar^{-1} \e^{-\xi/\xibar} = \xibar^{-1} - \xibar^{-2} \xi + O(\xi^2).
\end{equation}
The asymptotics of $k(\vecw',\xi,\vecw)$ for large $\xi$ is more complicated to state, see \cite{partIV}. We will here focus on tail asymptotics for the distribution of free path lengths \cite{partIV}. For any $\xi>0$, we have
\begin{equation}
\Psi_0(\xi) = \frac{1}{\xibar\zeta(d)} +O(\xi) ,
\end{equation}
and for $\xi\to\infty$
\begin{equation}\label{large}
\Psi_0(\xi) = \frac{A_d}{\xi^3} +O\bigl(\xi^{-3-\frac 2d}\bigr)
\begin{cases}
1&\text{if }\:d=2 
\\
\log\xi&\text{if }\:d=3
\\
1&\text{if }\:d\geq4 
\end{cases}
\end{equation}
with the constant
\begin{equation}
A_d=\frac{2^{2-d}}{d(d+1)\zeta(d)}.
\end{equation}
These asymptotics sharpen earlier upper and lower bounds by Bourgain, Golse and Wennberg \cite{Bourgain98,Golse00}.
Note that \eqref{large} implies that the density $\Psi_0(\xi)$ has no second moment. In dimension $d=2$ there is an explicit formula for $\Psi_0(\xi)$ conjectured by Dahlqvist \cite{Dahlqvist97}, and proved by Boca and Zaharescu \cite{Boca07}. This formula of course also follows directly from the expression for the transition kernel \eqref{twodidi}, cf.~\cite{partIII}.

\section{Several lattices} \label{sec:union} 

The previous two examples, random and periodic, could be analysed without the need to introduce colour. We will now describe a first example where the extension of $\Xi(t)$ to a Markov process $\widetilde\Xi(t)$ (as outlined in Section \ref{sec:generalized}) requires finitely many colours. 

We consider a scattering configuration given by the union of $N$ distinct affine Euclidean lattices,  
\begin{equation}
\scrP = \bigcup_{i=1}^N \scrL_i
\end{equation}
where each $\scrL_i$ has covolume $\nbar_i^{-1}$.
We will assume that the lattices are {\em pairwise incommensurable} in the sense that for any $i\neq j$, $c>0$ and $\veca\in\RR^d$, the intersection
$\scrL_i\cap(c\scrL_j+\veca)$ is contained in some affine linear subspace of dimension strictly less than $d$.\footnote{This condition is not essential in the proof of convergence, but ensures that the limit distributions have a particularly simple form. The case when all $N$ lattices are commensurable is a special case of the setting discussed in Section \ref{sec:quasicrystals}.} This ensures in particular that the density of $\scrP$ is $\nbar=\nbar_1+\ldots+\nbar_N$. 
As before, we stipulate without loss of generality that $\nbar=1$. 

To describe the random point processes and corresponding collision kernels, we require, in addition to a random lattice $\Theta_{\text{\rm lattice}}$ in the previous section, the notion of a random {\em affine} lattice. This is defined as 
$\Theta_{\text{\rm affine}}=(\ZZ^d+\vecalf)M$ where $\vecalf$ is a random variable uniformly distributed in $\TT^d=\ZZ^d\backslash\RR^d$ and $M$ is distributed with respect to Haar measure $\mu$ on $\GamG$ as before. Note that $\Theta_{\text{\rm affine}}$ is well defined, since $\TT^d$ and the Lebesgue measure on $\TT^d$ are invariant under the natural $\Gamma$ action (by right multiplication). We denote by $\Theta_{\text{\rm affine}}^{(1)},\ldots,\Theta_{\text{\rm affine}}^{(N)}$ independent copies of $\Theta_{\text{\rm affine}}$, which are furthermore independent of $\Theta_{\text{\rm lattice}}$. 

For $\vecy\in\scrL_j$ for some $j$, and $\vecy\notin\scrL_i$ for all $i\neq j$, we define the point process $\Theta_{\text{\rm union}}(\vecy)$ by
\begin{equation}
\Theta_{\text{\rm union}}(\vecy) = \nbar_j^{-1/d} \Theta_{\text{\rm lattice}} \cup \bigcup_{i\neq j} \big( \nbar_i^{-1/d} \Theta_{\text{\rm affine}}^{(i)} \big).
\end{equation}
In the following theorem, we say $\vecy\in\scrP$ is {\em generic}, if $\vecy\in\scrL_j$ is not rationally related to the other lattices $\scrL_i$ ($i\neq j$) in a sense made precise in \cite{union} (see the discussion after \cite[Thm.~1]{union}). The set of non-generic $\vecy$ in $\scrP$ is contained in a finite union of affine subspaces of dimension $<d$, and hence has zero relative density.

\begin{thm}[\cite{union}]\label{thm13.111}
Let $\lambda$ be an absolutely continuous probability measure on $\UB$, $\scrA_1,\ldots,\scrA_k\subset\RR^d$ bounded with boundary of measure zero and $n_1,\ldots,n_k\in\ZZ_{\geq 0}$. Then, for generic $\vecy\in\scrP$,
\begin{equation}\label{statem222}
\lim_{r\to 0}\PP(\#(\widetilde\Theta_r(\vecy)\cap\scrA_i) = n_i\,\forall i)  = \PP\big( \#((\Theta_{\text{\rm union}}(\vecy)-(0,\vech))\cap \scrA_i)=n_i \,\forall i  \big) .
\end{equation}
\end{thm}

The current setting requires $N$ colours. In the notation of Section \ref{sec:beyond}, we set $\Sigma=\{1,\ldots,N\}$,
$\iota(\vecy)=i$ if $\vecy\in\scrL_i$,
and define $\mm$ as the probability measure on $\Sigma$ so that $\mm(\{i\})=\nbar_i$. We prove in \cite{union} that the probability of emerging from a generic (as defined above) scatterer with a given colour $j'$ and random exit parameter $\vecw'$ (distributed according to a fixed, absolutely continuous Borel probability measure $\lambda$ on $\UB$), and hitting the next scatterer at time $T_n\in\;]\xi,\xi+d\xi[$ with colour $j$ and impact parameter $\vecw\in B \subset\UB$ converges in the Boltzmann-Grad limit to \eqref{sat9}. If the lattices are incommensurable as assumed above, the transition kernel in \eqref{sat9} is given by
\begin{equation}
k((\vecw',j),\xi,(\vecw,j)) =
k^{(1)}(\vecw',\nbar_j\xi,\vecw) \prod_{\substack{i=1 \\ i\neq j}}^N
\int_{\nbar_i \xi}^\infty \Psi(\xi')\, d\xi' ,
\end{equation}
and for $j'\neq j$,
\begin{equation}
k((\vecw',j'),\xi,(\vecw,j)) = \xibar\,
K^{(1)}(\nbar_{j'}\xi,\vecw')\, K^{(1)}(\nbar_j\xi,\vecw) \prod_{\substack{i=1 \\ i\neq j'\!\!,j}}^N
\int_{\nbar_i \xi}^\infty \Psi(\xi')\, d\xi' ,
\end{equation}
where $k^{(1)}(\vecw',\xi,\vecw)$ is the transition kernel for a single lattice in \eqref{1st}, $K^{(1)}(\xi,\vecw)$ the corresponding integrated kernel in \eqref{bigK} for a single lattice, and 
\begin{equation}
\Psi(\xi) := \frac{1}{v_{d-1}} \int_{\scrB_1^{d-1}} K^{(1)}(\xi,\vecw)\, d\vecw .
\end{equation}

The above formulas and \eqref{large} imply the following tail estimate for the distribution of free path lengths:
\begin{equation}
\Psi_0(\xi) = 
\frac{N(N+1) A_d^N \overline\sigma^{N-1}}{2^N \nbar_1\cdots\nbar_N} \; \xi^{-(N+2)} \times
\begin{cases}
\bigl(1+O(\xi^{-1})\bigr) &\text{if }\:d=2 
\\
\bigl(1+O(\xi^{-\frac23}\log\xi)\bigr) &\text{if }\:d=3
\\
\bigl(1+O(\xi^{-\frac2d})\bigr) &\text{if }\:d\geq4 .
\end{cases}
\end{equation}

The proof of the above results follows the same strategy as in the single-lattice case studied in Section \ref{sec:periodic}. The principal difference is that the equidistribution in the space of lattices stated in Theorem \ref{thm13.2} has to be generalised to the equidistribution in products:
Consider the subgroup $\widehat\Gamma=\Gamma_1\times\cdots\times\Gamma_N$ in $\SLdR^N$,
where each $\Gamma_i$ is a lattice in $\SLdR$.
We denote by $\mu_{\widehat\Gamma}$ the unique $\SLdR^N$ invariant probability measure on $\widehat\Gamma\bs \SLdR^N$, and by $\varphi$ the diagonal embedding of $\SLdR$ in $\SLdR^N$, i.e. $\varphi(M)=(M,\ldots,M)$. 
Recall that two lattices $\Gamma$ and $\Gamma'$ in $\SL(d,\RR)$ are said to be \textit{commensurable} if
their intersection $\Gamma\cap\Gamma'$ is also a lattice;
otherwise $\Gamma$ and $\Gamma'$ are \textit{incommensurable}.

\begin{thm}[\cite{union}]\label{equiThm2}
Let $\Gamma_1,\ldots,\Gamma_N\in\SLdR$ be pairwise incommensurable lattices, and $M\in\SLdR$. Let $\lambda$ be a Borel probability measure on $\UB$, absolutely continuous with respect to Lebesgue measure, and let $f:\UB\times \widehat\Gamma\bs \SLdR^N\to\RR$ be bounded continuous. Then
\begin{multline}\label{limiteq}
	\lim_{r\to 0} \int_{\UB} f\big(\vecw,\varphi(M S(\vecw)D(r))\big)\, d\lambda(\vecw)  \\
	= \int_{\UB\times\widehat\Gamma\bs \SLdR^N} f(\vecw,g) \, d\lambda(\vecw) \, d\mu_{\widehat\Gamma}(g) .
\end{multline}
\end{thm}

The key ingredient in the proof of this statement is Ratner's measure classification theorem \cite{Ratner91a} via a theorem of Shah on the equidistribution of translates of unipotent orbits \cite[Thm.~1.4]{Shah96}. Theorem \ref{thm13.2} corresponds to the special case $N=1$. For $N=2$ the proof is simpler than for $N\geq 3$, see \cite{farey}. Theorem \ref{equiThm2} is in fact an oversimplification---the proof of convergence to the transition kernel $k(\omega',\xi,\omega)$ in fact requires a variant of Theorem \ref{equiThm2} for products of spaces of affine lattices, cf.~\cite[Thm.~10]{union}.

The paper \cite{union} proves the convergence to $k(\omega',\xi,\omega)$ for a random exit parameter with fixed probability measure $\lambda$. What is still missing is a proof of the analogue of Theorem \ref{BBS} (for random scatterer configurations $\scrP$) or Theorem \ref{thm3} (where $\scrP$ is a single lattice). It is likely that the proof will follow the same line of arguments as in the periodic setting \cite{partII}.

\section{Quasicrystals}\label{sec:quasicrystals}

The third class of examples for  scattering configurations $\scrP$ that lead to a generalised Boltzmann equation---and the second that requires colour---are {\em quasicrystals}. We restrict our attention to quasicrystals constructed by the {\em cut-and-project method}, following closely the presentation in \cite{quasi}. Examples include many classic quasicrystals (such as the vertex set of a Penrose tiling) as well as locally finite periodic point sets. In contrast to the previous section, cut-and-project scatterer configurations generally require a continuous spectrum of colours. 

A cut-and-project set $\scrP\subset\RR^d$ is defined as follows, cf.~\cite{Baake}. For $m\geq 0$, $n=d+m$, let
\begin{equation}
\pi : \RR^{n} \to \RR^d ,
\qquad
\pi_\intl : \RR^{n} \to \RR^m 
\end{equation}
be the orthogonal projections of $\RR^{n}=\RR^d\times\RR^m$ onto the first and second factor, respectively. $\RR^d$ will be called the {\em physical space}, and $\RR^m$ the {\em internal space}. Let $\scrL\subset\RR^{n}$ be a lattice of full rank. The closure
\begin{equation}
\scrA := \overline{\pi_\intl(\scrL)} \subset\RR^m
\end{equation}
is an abelian subgroup, and we denote by $\scrA^0$ the connected component of $\scrA$ containing $0$. $\scrA^0$ is a linear subspace of $\RR^m$ of dimension $m_1$. We find vectors $\veca_1,\ldots,\veca_{m_2}$ ($m=m_1+m_2$) so that
\begin{equation}
\scrA = \scrA^0 \oplus \ZZ\pi(\veca_1)\oplus \ldots \oplus \ZZ\pi(\veca_{m_2}) .
\end{equation}
The Haar measure of $\scrA$ is denoted by $\mu_\scrA$ and normalised so that $\mu_\scrA\big|_{\scrA^0}$ is the standard Lebesgue measure on $\scrA^0$.
For $\scrV:=\RR^d\times \scrA^0$, we note that $\scrL\cap\scrV$ is a full rank lattice in $\scrV$. For $\scrW\subset\scrA$ with non-empty interior, we call 
\begin{equation}
\scrP=\scrP(\scrW,\scrL)=\{ \pi(\vecell) : \vecell\in\scrL, \; \pi_\intl(\vecell)\in\scrW \}
\end{equation}
a {\em cut-and-project set}. $\scrW$ is called the {\em window set}. If the boundary of the window set has $\mu_\scrA$-measure zero, we say $\scrP(\scrW,\scrL)$ is {\em regular}. We will furthermore assume that $\scrW$ and $\scrL$ are chosen so that the map
\begin{equation}
\pi_\scrW : \{ \vecell\in\scrL : \pi_\intl(\vecell)\in\scrW\} \to \scrP
\end{equation}
is bijective. This is to avoid coincidences in $\scrP$. It follows from Weyl equidistribution that such $\scrP$ have density
\begin{equation}
\nbar = \frac{\mu_\scrA(\scrW)}{\vol_{\RR^d}(\scrV/(\scrL\cap\scrV))} .
\end{equation}
Furthermore, for $\vecy\in\scrP$ there is $\vecell\in\scrL$ such that $\vecell=\pi(\vecy)$ and
\begin{equation}\label{eq16.8}
\scrP(\scrW,\scrL) -\vecy = \scrP(\scrW-\vecy_\intl,\scrL) , \qquad \vecy_\intl:=\pi_\intl(\vecell).
\end{equation}
This suggests to define the colour chart $\iota: \scrP\to \Sigma:=\scrW$ with $\iota(\vecy) = \vecy_\intl$.
The aim is now to describe the ``closure'' (in a suitable sense) of the orbit of $\scrP$ under the $\SLdR$-action and construct a probability measure on it. This will yield, as we shall see, our limit random process $\Theta(\vecy)$ in \eqref{limitTheta}. 

Set $G=\SL(n,\RR)$, $\Gamma=\SL(n,\ZZ)$ and define the embedding (for any $g\in G$) 
\begin{equation}
\varphi_g  : \SLdR  \hookrightarrow G ,\qquad 
A  \mapsto g \begin{pmatrix} A & 0_{d\times m} \\ 0_{m\times d} & 1_m \end{pmatrix} g^{-1} .
\end{equation}
Since $\SLdR$ is generated by unipotent subgroups, Ratner's theorems \cite{Ratner91a,Ratner91b} imply that there is a (unique) closed connected subgroup $H_g\leq G$ such that:
\begin{enumerate}[(i)]
\item $\Gamma\cap H_g$ is a lattice in $H_g$.
\item $\varphi_g(\SLdR)\subset H_g$.
\item The closure of $\Gamma\backslash\Gamma\varphi_g(\SLdR)$ is $\Gamma\backslash\Gamma H_g$.
\end{enumerate}
We will call $H_g$ a {\em Ratner subgroup.} We denote the unique $H_g$-invariant probability measure on $\Gamma\backslash\Gamma H_g$ by $\mu_{H_g}=\mu_g$.
Note that $\Gamma\backslash\Gamma H_g$ is isomorphic to the homogeneous space $(\Gamma\cap H_g)\backslash H_g$.

Pick $g\in G$, $\delta>0$ such that $\scrL=\delta^{1/n} \ZZ^n g$. Then one can show \cite[Prop.~3.5]{quasi} that $\pi_\intl(\delta^{1/n} \ZZ^n h g) \subset \scrA$ for all $h\in H_g$, and $\pi_\intl(\delta^{1/n} \ZZ^n h g) = \scrA$ for $\mu_g$-almost all $h\in H_g$.
The image of the map
\begin{equation}
\Gamma\backslash \Gamma H_g  \to \{\text{point sets in $\RR^d$}\} ,\qquad 
h  \mapsto \scrP(\scrW-\vecy_\intl, \delta^{1/n} \ZZ^n h g)
\end{equation}
defines a space of cut-and-project sets, and the push-forward of $\mu_g$ equips it with a probability measure. We have thus defined a random point process $\Theta_{\text{\rm quasi}}(\vecy)$ in $\RR^d$, which is $\SLdR$ invariant, and whose typical realisation is a cut-and-project set with window $\scrW-\vecy_\intl$ and internal space $\scrA$. This process is precisely the limit process we are looking for:

\begin{thm}[\cite{quasi}]\label{thm18.2}
Let $\lambda$ be an absolutely continuous probability measure on $\UB$, $\scrA_1,\ldots,\scrA_k\subset\RR^d$ bounded with boundary of measure zero and $n_1,\ldots,n_k\in\ZZ_{\geq 0}$. Then, for every $\vecy\in\scrP(\scrW,\scrL)$,
\begin{equation}\label{statem333}
\lim_{r\to 0} \PP(\#(\widetilde\Theta_r(\vecy)\cap\scrA_i) = n_i \,\forall i)  = \PP\big( \#((\Theta_{\text{\rm quasi}}(\vecy)-(0,\vech))\cap \scrA_i)=n_i \, \forall i \big) .
\end{equation}
\end{thm}

This statement is (as in previous sections) a consequence of equidistribution.
The following equidistribution theorems generalise Theorem \ref{thm13.2} stated earlier, and are used in the proof of Theorem \ref{thm18.2}. As in the case of Theorem \ref{equiThm2}, they are a consequence of Ratner's measure classification theorems \cite{Ratner91a}, and in particular follow from a theorem of Shah \cite[Thm.~1.4]{Shah96} on the equidistribution of translates of unipotent orbits. 

\begin{thm}[\cite{quasi}]\label{thm18.1}
Fix $g\in G$, $M\in\SLdR$.
For any bounded continuous $f: \UB \times \Gamma\backslash \Gamma H_g \to \RR$ and any absolutely continuous probability measure $\lambda$ on $\UB$, 
 \begin{equation}\label{eq18.2i}
\lim_{r\to 0} \int_{\UB} f(\vecw, \varphi_g(M S(\vecw) D(r)))\, d\lambda(\vecw) =
\int_{\UB} \int_{\Gamma\backslash \Gamma H_g} f(\vecw,h) \,d\mu_g(h)\, d\lambda(\vecw).
\end{equation}
\end{thm}

What are the subgroups $H_g$ that can arise in the above construction? For almost every lattice $\scrL$ in the space of lattices, we have $H_g=G$. Furthermore, if $m<d$, then for {\em every} $\scrL$ with the property that $\pi|_\scrL$ is injective, we have $H_g=G$ \cite[Prop.~2.1]{quasi}. A interesting class of examples when $m\geq d$ and $H_g\neq G$ are cut-and-project sets constructed from algebraic number fields. The  Penrose tilings fall into this class. Let us briefly sketch how such quasicrystals can be obtained as cut-and-project sets. Let $K$ be a totally real number field of degree $N\geq 2$ over $\QQ$, $\fO_K$ the ring of integers of $K$, and $\pi_1=\id$, $\pi_2,\ldots,\pi_N$ the distinct embeddings $K\hookrightarrow\RR$.
We also use $\pi_i$ to denote the component-wise embeddings
\begin{equation}
\pi_i :  K^d \hookrightarrow \RR^d ,\qquad \vecx \mapsto (\pi_i(x_1),\ldots ,\pi_i(x_d)) ,
\end{equation}
and similarly for the entry-wise embeddings of $d\times d$ matrices,
\begin{equation}
\pi_i :  \M_d(K) \hookrightarrow \M_d(\RR). 
\end{equation}
Now consider the lattice 
\begin{equation}
\scrL=\{ (\vecx, \pi_2(\vecx),\ldots,\pi_N(\vecx) ) :\vecx\in\fO_K^d \}
\end{equation}
in $\RR^{Nd}$. This is a lattice of full rank. The dimension of the internal space is $m=(N-1)d$. It is a fact of ``basic'' number theory \cite{Weil74} that $\scrA:=\overline{\pi_\intl(\scrL)} =\RR^m$, so that $\scrV=\RR^{Nd}$. Choose $g\in G$ and $\delta>0$ so that $\scrL=\delta^{1/Nd} \ZZ^{Nd} g$. Then \cite[Sect.~2.2.1.]{quasi} shows that 
\begin{equation}
H_g = g \SLdR^N g^{-1},\qquad
\Gamma\cap H_g = g \SL(d,\fO_K) g^{-1} ,
\end{equation}
where $\SL(d,\fO_K)$ is a Hilbert modular group.

A further example of a cut-and-project set is to take the union of finite translates of a given cut-and-project set. This is explained in \cite[Sect.~2.3]{quasi}. Let us here discuss the special case of periodic Delone sets, i.e., the union finite translates of a given lattice $\scrL_0$ of full rank in $\RR^d$. An example of such a set is the honeycomb lattice, which in the context of the Boltzmann-Grad limit of the Lorentz gas was recently studied by Boca et al.~\cite{Boca09,Boca10} with different techniques. The scatterer configuration $\scrP$ we are now interested in is the union of $m$ copies of the same lattice $\scrL_0$ translated by $\vect_1,\ldots,\vect_m\in\RR^d$,
\begin{equation}
\scrP = \bigcup_{j=1}^m (\vect_j + \scrL_0) .
\end{equation}
We assume that the $\vect_j$ are chosen in such a way that the above union is disjoint. Let us now show that $\scrP$ can be realised as a cut-and-project set $\scrP(\scrL,\scrW)$. Let
\begin{equation}
\scrL = (\scrL_0 \times \{\vecnull\}) + \sum_{j=1}^m \ZZ\, (\vect_j,\vece_j) \subset\RR^{n},
\end{equation}
where $\vecnull\in\RR^m$ and $\vece_1,\ldots,\vece_m$ are the standard basis vectors in $\RR^m$. The set $\scrL$ is evidently a lattice of full rank in $\RR^{n}$. Note that 
\begin{equation}
\pi_\intl(\scrL) = \sum_{j=1}^m \ZZ\, \vece_j = \ZZ^m,
\end{equation}
and therefore the closure of this set is $\scrA=\ZZ^m$ with connected component $\scrA^0=\{\vecnull\}$.
It follows that for the window set
\begin{equation}
\scrW= \bigcup_{j=1}^m \{ \vece_j \} \subset\scrA
\end{equation}
we indeed have 
\begin{equation}
\scrP(\scrL,\scrW) = \bigcup_{j=1}^m (\vect_j + \scrL_0) .
\end{equation}

Let us now determine $H_g$ in this setting. Take $g_0\in\SLdR$ so that $\scrL_0=\nbar_0^{-1/d}\ZZ^d g_0$, where $\nbar_0$ is the density of $\scrL_0$. Set
\begin{equation}
T = \begin{pmatrix} \vect_1 \\ \vdots \\ \vect_m \end{pmatrix} \in \M_{m\times d}(\RR).
\end{equation}
We then have $\scrL=\nbar_0^{-1/n}\ZZ^n g$,
for
\begin{equation}
g=\nbar_0^{1/n}\begin{pmatrix} \nbar_0^{-1/d} g_0 &  0 \\ T & 1_m \end{pmatrix}
\in\SL(n,\RR).
\end{equation}
Suppose $\veca_1,\ldots,\veca_d$ is a basis of $\scrL_0$ so that the vectors $\veca_1,\ldots,\veca_d$, $\vect_1,\ldots,\vect_m$ are linearly independent over $\QQ$. Then 
\begin{equation}
H_g = \bigg\{ \begin{pmatrix} h & 0 \\ u & 1_m \end{pmatrix} : h\in\SLdR,\; u\in\M_{m\times d}(\RR) \bigg\} .
\end{equation}
The Ratner subgroups that appear in the case of rational translates $\vect_j$ are discussed in \cite[Sect.~2.3.1]{quasi}.

Theorem \ref{thm18.2} gives a complete description of the limit processes $\Theta(\vecy)$ that may arise in the case of cut-and-project sets (as defined above). This answers in particular a question on the distribution of free path lengths raised by Wennberg \cite{Wennberg12}, see \cite{quasi} for details. We do not have a comprehensive solution to the remaining  {\em ``Does the limit \eqref{Lorentz-limit} exist?''}
and {\em ``What is the transition kernel $k(\omega',\xi,\omega)$?''} yet, but plan to address these in a forthcoming paper \cite{quasikinetic}.

\section{Superdiffusion}\label{sec:superdiffusion}

One of the central challenges in non-equilibrium statistical mechanics is to establish whether the dynamics of a test particle converges, in the limit of large times and after a suitable rescaling of length units, to Brownian motion. The first important step in the proof of such an invariance principle is the central limit theorem for the displacement $\vecQ(t)-\vecQ_0$, suitably normalised by a factor $\sigma(t)$. If $\sigma(t)\asymp \sqrt t$, we say the dynamics is {\em diffusive}. If $\sigma(t)/\sqrt t\to 0$ or $\sigma(t)/\sqrt t\to\infty$ as $t\to\infty$, the dynamics is called {\em subdiffusive} or {\em superdiffusive}, respectively. In the case of fixed scatterer radius $r$, most results are restricted to the periodic setting and dimension $d=2$, recall Section \ref{sec:introduction}. In the case of the Boltzmann-Grad limit with a random scatterer configuration, we have a central limit theorem with standard $\sqrt t$ normalisation:

\begin{thm}\label{thm:main00}
Let $\vecQ^\BG(t)$ denote the position variable of the random flight process $\Xi(t)$ for a Poisson scatterer configuration (cf.~Section \ref{sec:Poisson}). Then there exists a constant $\sigma_d>0$ such that, for any bounded continuous $f:\RR^d\to\RR$ and any\footnote{Because we have already passed to the Boltzmann-Grad limit, we may here consider the random process $\vecQ^\BG(t)$ either with {\em fixed} initial data (as stated) or with random initial data distributed according to $\Lambda$ (as assumed in all previous sections).\label{foot19}} $(\vecQ_0,\vecV_0)\in\T^1(\RR^d)$,
\begin{equation}\label{eq:main00}
\lim_{t\to\infty}\Exp f\bigg(\frac{\vecQ^\BG(t) -\vecQ_0}{\sigma_d \sqrt{t}}  \bigg) =
\frac{1}{(2\pi)^{d/2}} \int_{\RR^d} f(\vecx)\, \e^{-\frac12 \|\vecx\|^2} d\vecx .
\end{equation}
\end{thm}

This theorem follows from standard techniques in the theory of Markov processes \cite{Papanicolaou75}, as pointed out by Spohn \cite{Spohn78}. On the other hand, the Boltzmann-Grad limit of a periodic Lorentz gas satisfies a superdiffusive central limit theorem with $\sqrt{t\log t}$ normalisation:

\begin{thm}[Marklof and T\'oth, 2014 \cite{super}]\label{thm:main1}
Let $\vecQ^\BG(t)$ denote the position variable of the random flight process $\Xi(t)$ for a periodic scatterer configuration (cf.~Section \ref{sec:periodic}). Then, for any bounded continuous $f:\RR^d\to\RR$ and any\footnote{Cf.~footnote \ref{foot19}.} $(\vecQ_0,\vecV_0)\in\T^1(\RR^d)$,
\begin{equation}\label{eq:main1}
\lim_{t\to\infty}\Exp f\bigg(\frac{\vecQ^\BG(t) -\vecQ_0}{\Sigma_d \sqrt{t\log t}}  \bigg) =
\frac{1}{(2\pi)^{d/2}} \int_{\RR^d} f(\vecx)\, \e^{-\frac12 \|\vecx\|^2} d\vecx 
\end{equation}
with $\Sigma_d^2:=\frac{A_d}{2d\xibar}$.
\end{thm}

Recall that $A_d$ is the constant in the tail asymptotics of the free path lengths \eqref{large}. This means in particular that $\Sigma_d$ is independent of the choice of scattering map (within the admissible class). Although the superdiffusive scaling is intimately related to the fact that the second moment of the distribution of free path lengths diverges, the proof of Theorem \ref{thm:main1} requires further information on the transition kernel $k(\omega',\xi,\omega)$. The main ingredients of our proof are (a) exponential decay of correlations in the sequence of random variables $(\eta_n,\vecV_n)_{n\in\NN}$ and (b) the Lindeberg central limit theorem for the independent random variables $(\xi_n | \ueta)_{n\in\NN}$ conditioned on $\ueta=(\eta_n)_{n\in\NN}$. For full details, see \cite{super}.

\section*{Acknowledgment}

Much of the material presented in this paper is based on joint work with Andreas Str\"ombergsson and B\'alint T\'oth. I would like to thank Andreas and B\'alint for the fruitful collaboration. I am grateful to Daniel El-Baz, Jory Griffin, Andreas Str\"ombergsson, B\'alint T\'oth, Jim Tseng and Ilya Vinogradov for their comments on the first draft of this paper.

\end{document}